\documentclass{amsart}
\usepackage{amsmath,amsthm,amssymb,amsfonts,amscd, array}
\newtheorem{Theorem}{Theorem}[section]
\newtheorem{defn}[Theorem]{Definition}
\newtheorem{cor}[Theorem]{Corollary}
\newtheorem{prp}[Theorem]{Proposition}
\newtheorem{lemma}[Theorem]{Lemma}
\newtheorem{Example}[Theorem]{Example}
\newtheorem{rmk}[Theorem]{Remark}
\title[On $\phi$-$\delta$-$S$-primary ideals of Commutative rings]
{On $\phi$-$\delta$-$S$-primary ideals of Commutative rings}
\author[Ameer Jaber]{}
 \email{ameerj@hu.edu.jo}
 \medskip

\begin{document}

\maketitle

\centerline{\scshape  Ameer Jaber}
\medskip

{\footnotesize \centerline{ Department of Mathematics }
\centerline{ College of Science }
\centerline{ The Hashemite University} \centerline{ Zarqa 13115, Jordan
} }
\begin{abstract}
Let $R$ be a commutative ring with unity $(1\not=0)$ and let $\mathfrak{J}(R)$ be the set of all ideals of $R$. Let $\phi:\mathfrak{J}(R)\rightarrow\mathfrak{J}(R)\cup\{\emptyset\}$ be a reduction function of ideals of $R$ and let $\delta:\mathfrak{J}(R)\rightarrow\mathfrak{J}(R)$ be an expansion function of ideals of $R$. We recall that a proper ideal $I$ of $R$ is called a $\phi$-$\delta$-primary ideal of $R$ if whenever $a,b\in R$ and $ab\in I-\phi(I)$, then $a\in I$ or $b\in\delta(I)$. In this paper, we introduce a new class of ideals that is a generalization to the class of $\phi$-$\delta$-primary ideals. Let $S$ be a multiplicative subset of $R$ such that $1\in S$ and let $I$ be a proper ideal of $R$ with $S\cap I=\emptyset$, then $I$ is called a $\phi$-$\delta$-$S$-primary ideal of $R$ associated to $s\in S$ if whenever $a,b\in R$ and $ab\in I-\phi(I)$, then $sa\in I$ or $sb\in\delta(I)$. In this paper, we have presented a range of different examples, properties, characterizations of this new class of ideals.
\end{abstract}
\section{Introduction}
Throughout this paper all rings are commutative with $(1\not=0).$ Let $\mathfrak{J}(R)$ be the set of all ideals of $R$. A. Hamed and A. Malek ~\cite{HM} introduced the concept of $S$-prime ideals which is generalization of prime ideals, where $S$ is a is a multiplicative subset of $R$ such that $1\in S.$ Recall from ~\cite{HM} that a proper ideal $I$ of $R$ with $I\cap S=\emptyset$ is said to be $S$-prime if there exists $s\in S$ such that for all $a,b\in R$ with $ab\in I$ implies that $sa\in I$ or $sb\in I$. In a recent study, F. Almahdi, E. Bouba, and M. Tamekkante ~\cite{FBT} introduced the concept of weakly $S$-prime ideals, which is also a generalization of $S$-prime ideals, prime ideals and weakly prime ideals, where $S$ is a is a multiplicative subset of $R$ such that $1\in S.$ Recall from ~\cite{FBT} that a proper ideal $I$ of $R$ with $I\cap S=\emptyset$ is said to be weakly $S$-prime if there exists $s\in S$ such that for all $a,b\in R$ with $0\not=ab\in I$ implies that $sa\in I$ or $sb\in I$. D. Zhao ~\cite{DZ} introduced the concept of expansion function of ideals of $R$. Let $\delta$ be an expansion function of ideals of $R$, recall from ~\cite{DZ} that a proper ideal $I$ of $R$ is said to be a $\delta$-primary ideal of $R$, if $a,b\in R$ with $ab\in I$, then $a\in I$ or $b\in\delta(I)$. This concept of $\delta$-primary ideals is a generalization of the concepts of prime ideals and primary ideals.
Let $\delta$ be an expansion function of ideals of $R$ and $\phi$ a reduction function of ideals of $R$. In ~\cite{AJ} the author introduced the concept of $\phi$-$\delta$-primary ideals and this concept is a generalization of the concept of $\delta$-primary ideals in ~\cite{DZ}. So from this point of view, $\phi$-$\delta$-primary ideals generalize the concepts of prime ideals, weakly prime ideals, almost prime ideals, primary ideals, weakly primary ideals and almost primary ideals. For more generalizations of prime ideals, we refer to \cite{ST1, ST2, ST3}. Our main purpose in this paper is to extend the concept of $\phi$-$\delta$-primary ideals of $R$ to the concept of $\phi$-$\delta$-$S$-primary ideals of $R$, where $S$ is a multiplicative subset of $R$ such that $1\in S.$ This means that the concept of $\phi$-$\delta$-$S$-primary ideals is a generalization of the concepts of prime ideals, weakly prime ideals, almost prime ideals, primary ideals, weakly primary ideals, almost primary ideals, $S$-prime ideals, weakly $S$-prime ideals, almost $S$-prime ideals, $S$-primary ideals, weakly $S$-primary ideals and almost $S$-primary ideals. However, the converse is not true in general. (See Example~\ref{2amr-21}, Example~\ref{2amr-22} and Example~\ref{2amr-23}).  Let $\phi,\delta$ be a reduction function and an expansion function of ideals of $R$, respectively, and let $S$ be a multiplicative subset of $R$ such that $1\in S.$ In this paper, we call a proper ideal $I$ of $R$, with $I\cap S=\emptyset$, a $\phi$-$\delta$-$S$-primary ideal of $R$ associated to some $s\in S$ if whenever $a,b\in R$ such that $ab\in I-\phi(I)$, then $sa\in I$ or $sb\in\delta(I)$. Among many results in this paper, it is shown (Proposition 2.17) that if $I$ is a $\phi$-$\delta$-$S$-primary ideal of $R$ associated to some $s\in S$ which is not $S$-primary, then $I^2\subseteq\phi(I)$. Theorem 2.20 proves that a proper ideal $I$ of $R$ is a $\phi$-$\delta$-$S$-primary ideal of $R$ associated to some $s\in S$ if and only if for each $a\not\in(\delta(I):s)$ we have either $(I:a)\subseteq(I:s)$ or $(I:a)=(\phi(I):a).$ Similarly, in Theorem 2.21, we prove that a proper ideal $I$ of $R$ is a $\phi$-$\delta$-$S$-primary ideal of $R$ associated to some $s\in S$ if and only if for each $a\not\in(I:s)$ we have either $(I:a)\subseteq(\delta(I):s)$ or $(I:a)=(\phi(I):a).$ In the case when $S$ satisfies the conditions $\phi(I)=(\phi(I):s)$ for some $s\in S$ and $(\phi(I):t)\subseteq(\phi(I):s)$ $\forall t\in S$, it is proved (Theorem 2.33) that $I$ is a $\phi$-$\delta$-$S$-primary ideal of $R$ associated to some $s\in S$ if and only if $(I:s)$ is a $\phi$-$\delta$-primary ideal of $R$ if and only if $S^{-1}I$ is a $\phi_S$-$\delta_S$-primary ideal of $S^{-1}R$ and $S^{-1}I\cap R=(I:s).$\\ In the next section, let $f:X\rightarrow Y$ be a nonzero $(\delta,\phi)-(\gamma,\psi)$-epimorphism. In Theorem 3.3, we prove that $f$ induces one-to-one correspondence between $\phi$-$\delta$-$S$-primary ideals of $X$ associated to some $s\in S$ consisting $\ker(f)$ and $\psi$-$\gamma$-$f(S)$-primary ideal of $Y$ associated to $f(s)\in f(S).$ Also, in Lemma 3.7, we prove that if $a,b\in X$, then $(a,b)$ is a $\phi$-$\delta$-$S$-twin zero of $I$, where $I$ is a $\phi$-$\delta$-$S$-primary ideals of $X$ associated to some $s\in S$ consisting $\ker(f)$, if and only if $(f(a),f(b))$ is a $\psi$-$\gamma$-$f(S)$-twin zero of $f(I).$\\ In the last section, we determine all $\phi$-$\delta$-$S$-primary ideals in direct product of rings and we prove some results concerning $\phi$-$\delta$-$S$-primary ideals in direct product of rings. (See, Theorems 4.1-4.4).
\section{Properties of $\phi$-$\delta$-$S$-Primary ideals}
\begin{defn}
Let $R$ be a commutative ring with unity ($1\not=0$), and let $\mathfrak{J}(R)$ be the set of all ideals of $R$.\\
(1) Recall from \cite{DZ} that a function $\delta : \mathfrak{J}(R)\rightarrow\mathfrak{J}(R)$ is called an expansion function of ideals of $R$ if whenever $I, J, K$ are ideals of $R$ with $J\subseteq I$, then $\delta(J)\subseteq\delta(I)$ and $K\subseteq\delta(K)$.\\
(2) Recall from \cite{AJ} that a function $\phi:\mathfrak{J}(R)\rightarrow\mathfrak{J}(R)$ is called a reduction function of ideals of $R$ if $\phi(I)\subseteq I$ for all ideals $I$ of $R$ and if whenever $P\subseteq Q$, where $P$ and $Q$ are ideals of $R$, then $\phi(P)\subseteq\phi(Q).$
\end{defn}
\begin{defn}
Let $R$ be a commutative ring with unity ($1\not=0$), and $S$ a multiplicative subset of $R$. Suppose $\delta$, $\phi$ are expansion and reduction functions of ideals of $R$, respectively.\\
 (1) A proper ideal $I$ of $R$ satisfying $I\cap S=\emptyset$ is said to be a $\delta$-$S$-primary ideal of $R$ associated to $s\in S,$ if whenever $ab\in I$, then $sa\in I$ or $sb\in\delta(I)$ for all $a,b\in R$.\\
 (2) A proper ideal $I$ of $R$ satisfying $I\cap S=\emptyset$ is said to be a $\phi$-$\delta$-$S$-primary ideal of $R$ associated to $s\in S,$ if whenever $ab\in I-\phi(I)$, then $sa\in I$ or $sb\in\delta(I)$ for all $a,b\in R$.
\end{defn}
Throughout this section, $R$ denotes a commutative ring with unity ($1\not=0$), $S$ denotes a multiplicative subset of $R$ such that $1\in S,$ $\delta,\gamma:\mathfrak{J}(R)\rightarrow\mathfrak{J}(R)$ denote expansion functions, and $\phi,\psi:\mathfrak{J}(R)\rightarrow\mathfrak{J}(R)$ denote reduction functions.\\ \\
In the following example, we recall from \cite{B} some examples of expansion functions of ideals of a given ring $R$.
 \begin{Example}\textnormal{$\ $\\
 (1) The identity function $\delta_0$, where $\delta_0(I)=I$ for any $I\in\mathfrak{J}(R)$, is an expansion function of ideals in $R$.\\
 (2) For each ideal $I$ of $R$ define $\delta_1(I)=\sqrt{I}$. Then $\delta_1$ is an expansion function of ideals in $R$.\\
 (3) Let $J$ be a proper ideal of $R$. If $\delta(I)=I+J$ for every ideal $I$ in $\mathfrak{J}(R)$, then $\delta$ is an expansion function of ideals in $R$.\\
 (4) Let $J$ be a proper ideal of $R$. If $\delta(I)=(I:J)$ for every ideal $I$ in $\mathfrak{J}(R)$, then $\delta$ is an expansion function of ideals in $R$.\\
 (5) Assume that $\delta_1$, $\delta_2$ are expansion functions of ideals of $R$. Let $\delta : \mathfrak{J}(R)\rightarrow\mathfrak{J}(R)$ such that $\delta(I) =\delta_1(I) + \delta_2(I)$. Then $\delta$ is an expansion function of ideals of $R$.\\
 (6) Assume that $\delta_1$, $\delta_2$ are expansion functions of ideals of $R$. Let $\delta : \mathfrak{J}(R)\rightarrow\mathfrak{J}(R)$ such that $\delta(I) =\delta_1(I)\cap\delta_2(I)$. Then $\delta$ is an expansion function of ideals of $R$.\\
 (7) Assume that $\delta_1$,...,$\delta_n$ are expansion functions of ideals of $R$. Let $\delta : \mathfrak{J}(R)\rightarrow\mathfrak{J}(R)$ such that $\delta(I) =\cap_{i=1}^n\delta_i(I)$
 then $\delta$ is also an expansion function of ideals of $R$.\\
 (8) Assume that $\delta_1$, $\delta_2$ are expansion functions of ideals of $R$. Let $\delta : \mathfrak{J}(R)\rightarrow\mathfrak{J}(R)$ such that $\delta(I) =\delta_1(\delta_2(I))$. Then $\delta$ is an expansion function of ideals of $R$.}
 \end{Example}
Recall that if $\psi_1, \psi_2 : \mathfrak{J}(R)\rightarrow\mathfrak{J}(R)\cup\{\emptyset\}$ are expansion (reduction) functions of ideals of $R$, we define
$\psi_1\leq\psi_2$ if $\psi_1(I)\subseteq\psi_2(I)$ for each $I\in\mathfrak{J}(R)$.\\
In the following example, we  recall from \cite{AB} some examples of reduction functions of ideals of a given ring $R$.
\begin{Example}\textnormal{$\ $\\
(1) The function $\phi_\emptyset$, where $\phi_\emptyset(I)=\emptyset$ for any $I\in\mathfrak{J}(R)$ is an ideal reduction.\\
(2) The function $\phi_0$, where $\phi_0(I)=\{0\}$ for any $I\in\mathfrak{J}(R)$ is an ideal reduction.\\
(3) The function $\phi_2$, where $\phi_2(I)=I^2$ for any $I\in\mathfrak{J}(R)$ is an ideal reduction.\\
(4) The function $\phi_n$, where $\phi_n(I)=I^n$ for any $I\in\mathfrak{J}(R)$ is an ideal reduction.\\
(5) The function $\phi_\omega$, where $\phi_\omega(I)=\cap_{n=1}^\infty I^n$ for any $I\in\mathfrak{J}(R)$ is an ideal reduction.\\
(6) The function $\phi_1$, where $\phi_1(I)=I$ for any $I\in\mathfrak{J}(R)$ is an ideal reduction.\\
Observe that $\phi_\emptyset\leq \phi_0\leq \phi_\omega\leq\cdots\leq \phi_{n+1}\leq \phi_n\leq\cdots\leq\phi_2\leq\phi_1$.}
\end{Example}
\begin{rmk}\textnormal{$\ $\\
(1) If $\delta\leq\gamma$. Then every $\phi$-$\delta$-$S$-primary ideal of $R$ is a $\phi$-$\gamma$-$S$-primary ideal. In particular, every $\phi$-$S$-prime ideal of $R$ is a $\phi$-$\delta$-$S$-primary ideal. However, the converse is not true in general.\\
(2) If $\phi\leq\psi$. Then every $\phi$-$\delta$-$S$-primary ideal of $R$ is a $\psi$-$\delta$-$S$-primary ideal. In particular, every $\delta$-$S$-primary ideal of $R$ is a $\phi$-$\delta$-$S$-primary ideal. However, the converse is not true in general.}
\end{rmk}
\begin{Example}\label{2amr-21}\textnormal{$\ $\\
(1) Set $R=\mathbb{Z}_{12}$, $I=4\mathbb{Z}_{12}.$ Then $\delta_1(I)=\sqrt{I}=2\mathbb{Z}_{12}$. Take $S=\{1\}$, $\phi=\phi_\emptyset.$ Then it is easy to check that $I$ is a $\delta_1$-$S$-primary ideal of $R.$ Moreover, $I$ is not an $S$-prime ideal, since $(2)(2)=4\in I$ but $2\not\in I.$\\
(2) Set $R=\mathbb{Z}_{12}$, $S=\{1,5\}$. Then $S$ is a multiplicative subset of $R$. Let $I=\{0\}.$ Then $\delta_1(I)=6\mathbb{Z}_{12},$ $\phi_2(I)=I^2=(0).$ So, $I$ is an almost-$\delta_1$-$S$-primary ideal of $R$ associated to $s=5.$ Moreover, $(3)(4)=0\in I$ but neither $(3)(5)=4\in\delta_1(I)$ nor $(4)(5)=8\in\delta_1(I)$. Thus, $I$ is not a $\delta_1$-$S$-primary ideal of $R$ associated to $s=5.$}
\end{Example}
\begin{prp}\textnormal{
Let $\{ J_i\ :\ i\in \vartriangle\}$ be a directed set of $\phi$-$\delta$-$S$-primary ideals of $R$ associated to $s\in S$. Then the ideal $J=\cup_{i\in \vartriangle}J_i$ is a $\phi$-$\delta$-$S$-primary ideal of $R$ associated to $s\in S$.}
\end{prp}
\textbf{Proof.}$\ $\\
Let $ab\in J-\phi(J)$, where $a, b\in R$. Suppose $sa\not\in J$. We want to show that $sb\in\delta(J)$. Since $ab\not\in\phi(J)$, we have $ab\not\in\phi(J_i)$ for all $i\in\vartriangle$. Let $t\in\vartriangle$ such that $ab\in J_t-\phi(J_t)$, then $sa\in J_t$ or $sb\in\delta(J_t)$, since $J_t$ is a $\phi$-$\delta$-$S$-primary ideal of $R$ associated to $s\in S$. Since $sa\not\in J$, we have $sa\not\in J_t$ which implies that $sb\in\delta(J_t)\subseteq\delta(J).$ Hence $J$ is a $\phi$-$\delta$-$S$-primary ideal of $R$ associated to $s\in S$.\hfill$\blacksquare$
\begin{prp}\textnormal{
Let $\{Q_i\ :\ i\in \vartriangle\}$ be a directed set of $\phi$-$\delta$-$S$-primary ideals of $R$ associated to $s\in S$. Suppose $\phi(Q_i)=\phi(Q_j)$ and $\delta(Q_i)=\delta(Q_j)$ for every $i,j\in\vartriangle$. If $\phi$, $\delta$ have the intersection property, then the ideal $J=\cap_{i\in \vartriangle}Q_i$ is a $\phi$-$\delta$-$S$-primary ideal of $R$ associated to $s\in S$.}
\end{prp}
\textbf{Proof.}$\ $\\
Let $t\in\vartriangle$, since $\phi(Q_i)=\phi(Q_t)$ and $\delta(Q_i)=\delta(Q_t)$ for every $i\in\vartriangle,$ and since $\phi$, $\delta$ have the intersection property, then $\phi(J)=\phi(Q_t)$ and $\delta(J)=\delta(Q_t).$ Let $ab\in J-\phi(J)$, where $a,b\in R$ such that $sb\not\in\delta(J)$. Then $ab\in Q_t-\phi(Q_t).$ Since $Q_t$ is a $\phi$-$\delta$-$S$-primary ideals of $R$ associated to $s\in S$, we conclude that $sa\in Q_t$ or $sb\in\delta(Q_t)$. Since $sb\not\in\delta(J)$, $sb\not\in\delta(Q_t)=\delta(J)$. Hence we conclude that $sa\in Q_t$ for each $t\in\vartriangle$ which implies that $sa\in J$. Thus, $J$ is a $\phi$-$\delta$-$S$-primary ideal of $R$ associated to $s\in S.$\hfill$\blacksquare$\\ \\
Obviously, every $\phi$-$\delta$-primary ideal $R$ is a $\phi$-$\delta$-$S$-primary ideal. In particular, every weakly-$\delta$-primary ($\delta$-primary) ideal of $R$ is a weakly-$\delta$-$S$-primary ($\delta$-$S$-primary). However, the next two examples show that the converse is not true in general.
\begin{Example}\label{2amr-22}\textnormal{
Let $R=\mathbb{Z}_{80}$, $I=20\mathbb{Z}_{80}$ and $S=\{1,5,25,45,65\}.$ Then $S$ is a multiplicative subset of $R$ such that $I\cap S=\emptyset.$ Let $\delta=\delta_1$ and $\phi=\phi_0$, then $\delta_1(I)=\sqrt{I}=10\mathbb{Z}_{80}$ and $\phi_0(I)=(0).$ Let $a, b\in R$ such that $0\not=ab\in I$, then $2/ab$ which implies that $2/a$ or $2/b.$ Thus, $5a\in\delta_1(I)$ or $5b\in\delta_1(I).$ Hence we conclude that $I$ is a weakly-$\delta_1$-$S$-primary ideal of $R$ associated to $s=5.$ Moreover, $0\not=(4)(5)=20\in I$ but neither $4\in\delta_1(I)$ nor $5\in\delta_1(I)$. Thus, $I$ is not a weakly-$\delta_1$-primary.}
\end{Example}
\begin{Example}\label{2amr-23}\textnormal{
Let $R=\mathbb{Z}[x]$, $I=<4x>=4x\mathbb{Z}[x].$ Let $\phi_2(I)=I^2=<16x^2>$ and $\delta_1(I)=\sqrt{I}=<2x>.$ Let $S=\{2^k:k\geq0\}$. Then $S$ is a multiplicative subset of $R$ such that $I\cap S=\emptyset.$ Moreover, $I$ is an almost-$\delta_1$-$S$-primary ideal of $R$ associated to $s=2\in S$, since if $f(x), g(x)\in R$ with $f(x)g(x)\in I-I^2$, then $x/f(x)$ or $x/g(x)$ which implies that $2f(x)\in\delta_1(I)$ or $2g(x)\in\delta_1(I)$. Since $4x\in I-I^2$ and neither $4\in\delta_1(I)$ nor $x\in\delta_1(I)$, then we get that $I$ is not an almost-$\delta_1$-primary ideal of $R.$}
\end{Example}
\begin{prp}\textnormal{
Let $I$ be a proper ideal of $R$ such that $I$ is a $\phi$-$S$-primary ideal of $R$ associated to $s\in S$ such that $\sqrt{\delta(I)}\subseteq\delta(\sqrt{I})$ and $\sqrt{\phi(I)}\subseteq\phi(\sqrt{I})$, then $\sqrt{I}$ is a $\phi$-$\delta$-$S$-primary ideal of $R$ associated to $s.$}
\end{prp}
\textbf{Proof.}$\ $\\
Let $a,b\in R$ such that $ab\in\sqrt{I}-\phi(\sqrt{I})$. Then $ab\in\sqrt{I}$ which implies that $a^nb^n\in I$ for some $n\geq1$. If $a^nb^n\in\phi(I)$, then $ab\in\sqrt{\phi(I)}\subseteq\phi(\sqrt{I})$, a contradiction. Thus, $a^nb^n\in I-\phi(I)$ which implies that $sa^n\in I$ or $sb^n\in\delta(I).$ Thus, $sa\in\sqrt{I}$ or $sb\in\sqrt{\delta(I)}\subseteq\delta(\sqrt{I})$. Hence, $\sqrt{I}$ is a $\phi$-$\delta$-$S$-primary ideal of $R$ associated to $s$.\hfill$\blacksquare$
\begin{cor}\textnormal{
Let $I$ be a proper ideal of $R$ such that $I$ is a $\phi$-$S$-primary ideal of $R$ associated to $s\in S.$ Suppose that $\sqrt{\phi(I)}\subseteq\phi(\sqrt{I})$. Then $\sqrt{I}$ is a $\phi$-$S$-prime ideal of $R$ associated to $s.$}
\end{cor}
\textbf{Proof.}$\ $\\
Let $\delta(J)=\sqrt{J}$ for ever ideal $J$ in $R$. Then, by the above proposition, if $I$ is a $\phi$-$S$-primary ideal of $R$ associated to $s$ then $\sqrt{I}$ is a $\phi$-$S$-prime ideal of $R$ associated to $s$.\hfill$\blacksquare$
\begin{prp}\textnormal{
Let $I$ be a proper ideal of $R$ such that $I$ is a $\phi$-$S$-primary ideal of $R$ associated to $s\in S.$ Suppose that $\sqrt{\phi(I)}\subseteq\phi(\sqrt{I})$ and $(\phi(\sqrt{I}):x)\subseteq(\phi(\sqrt{I}):s)$ for each $x\in S$. If $a\in R-(\sqrt{I}:s)$, then $S\cap(\sqrt{I}:a)=\emptyset.$}
\end{prp}
\textbf{Proof.}$\ $\\
It is easy to see that $\sqrt{I}\cap S=\emptyset$, since $I\cap S=\emptyset$. Also, by the above corollary, $\sqrt{I}$ is a $\phi$-$S$-prime ideal of $R$ associated to $s$. We show that $S\cap(\sqrt{I}:a)=\emptyset.$ Let $t\in S$ such that $ta\in\sqrt{I}.$ If $ta\in\phi(\sqrt{I}),$ then $a\in(\phi(\sqrt{I}):t)\subseteq(\phi(\sqrt{I}):s)$ which implies that $sa\in\phi(\sqrt{I})\subseteq\sqrt{I},$ a contradiction. Thus, $ta\in\sqrt{I}-\phi(\sqrt{I})$
implies that $sa\in\sqrt{I}$ or $st\in\sqrt{I},$ which is a contradiction again, since $a\not\in(\sqrt{I}:s)$ and $S\cap\sqrt{I}=\emptyset.$ Thus, $S\cap(\sqrt{I}:a)=\emptyset.$\hfill$\blacksquare$
\begin{cor}\textnormal{
Let $I$ be a proper ideal of $R$ such that $I$ is a $\phi$-$\delta$-$S$-primary ideal of $R$ associated to $s\in S$ with $\delta(I)\subseteq\sqrt{I}.$ Suppose $(\phi(\sqrt{I}):x)\subseteq(\phi(\sqrt{I}):s)$ for each $x\in S$ and $(\delta(I):s)=(\sqrt{I}:s)$. Then $(\delta(I):s)=(\delta(I):s^2)$ and if whenever $a\in R-(\delta(I):s)$, then $S\cap(\delta(I):a)=\emptyset.$}
\end{cor}
\textbf{Proof.}$\ $\\
Since $I$ is a $\phi$-$\delta$-$S$-primary ideal of $R$ associated to $s$ and $\delta(I)\subseteq\sqrt{I}$, it is easy to see that $I$ is a $\phi$-$S$-primary ideal of $R$ associated to $s$ and $(\sqrt{I}:s)=(\sqrt{I}:s^2).$ Since $\delta(I)\subseteq\sqrt{I}$ and $(\delta(I):s)=(\sqrt{I}:s),$ we have $(\delta(I):s)=(\delta(I):s^2)$. Moreover, if $a\in R-(\delta(I):s)$, then $sa\not\in\sqrt{I}.$ Thus, by the above proposition, $S\cap(\sqrt{I}:a)=\emptyset.$ Hence $S\cap(\delta(I):a)\subseteq S\cap(\sqrt{I}:a)=\emptyset$, since $\delta(I)\subseteq\sqrt{I}.$\hfill$\blacksquare$\\ \\
Recall that if $I,J,K$ are ideals of $R$ such that $K=I\cup J$, then $K=I$ or $K=J$.
\begin{Theorem}\label{2amr-1}\textnormal{
Let $I$ be a proper ideal of $R$ such that $I$ is a $\phi$-$\delta$-$S$-primary ideal of $R$ associated to $s\in S.$ If $a\in R-(\delta(I):s^2)$, then $(I:sa)=(I:s)$ or $(I:sa)=(\phi(I):sa).$}
\end{Theorem}
\textbf{Proof.}$\ $\\
 It is enough to show that $(I:sa)=(I:s)\cup(\phi(I):sa).$ It is easy to see that $(I:s)$ and $(\phi(I):sa)$ are subsets of $(I:sa)$. Let $r\in (I:sa)$, then $rsa\in I$. If $rsa\in\phi(I)$ then $r\in(\phi(I):sa).$ So we may assume that $rsa\not\in\phi(I)$. Thus, $rsa\in I-\phi(I)$ implies that $sr\in I$ since $s^2a\not\in\delta(I)$. So, $r\in(I:s)$. Thus, $(I:sa)=(I:s)\cup(\phi(I):sa).$ Hence $(I:sa)=(I:s)$ or $(I:sa)=(\phi(I):sa).$\hfill$\blacksquare$
\begin{cor}\textnormal{
Let $I$ be a proper ideal of $R$ such that $I$ is a $\phi$-$S$-primary ideal of $R$ associated to $s\in S.$ If $a\in R-(\sqrt{I}:s)$, then $(I:sa)=(I:s)$ or $(I:sa)=(\phi(I):sa).$}
\end{cor}
\textbf{Proof.}$\ $\\
Let $I$ be a proper ideal of $R$ such that $I$ is a $\phi$-$S$-primary ideal of $R$ associated to $s\in S.$ Then it is easy to see that $(\sqrt{I}:s)=(\sqrt{I}:s^2).$ Thus, if $a\in R-(\sqrt{I}:s)$, then $a\in R-(\sqrt{I}:s^2)$. Hence, by Theorem~\ref{2amr-1}, $(I:sa)=(I:s)$ or $(I:sa)=(\phi(I):sa).$\hfill$\blacksquare$
\begin{prp}\label{2amr-2}\textnormal{
Let $I$ be a proper ideal of $R$ such that $I$ is a $\phi$-$\delta$-$S$-primary ideal of $R$ associated to $s\in S.$ If $I$ is not a $\delta$-$S$-primary, then $I^2\subseteq\phi(I).$}
\end{prp}
\textbf{Proof.}$\ $\\
Suppose that $I^2\not\subseteq\phi(I).$ We claim that $I$ is a $\delta$-$S$-primary ideal of $R$ associated to $s$. Let $a,b\in R$ such that $ab\in I$. If $ab\in I-\phi(I)$, then $sa\in I$ or $sb\in\delta(I)$. Therefore we may assume that $ab\in\phi(I)$. Suppose that $aI\not\subseteq\phi(I)$, then there exists $p\in I$ such that $ap\not\in\phi(I)$. So, $a(p+b)\in I-\phi(I)$ implies that $sa\in I$ or $s(p+b)\in\delta(I),$ and since $sp\in I\subseteq\delta(I)$, we have $sb\in\delta(I)$. Similarly, if $bI\not\subseteq\phi(I)$, we obtain that $sa\in I$. Thus we may assume that $aI\subseteq\phi(I)$ and $bI\subseteq\phi(I).$ Since $I^2\not\subseteq\phi(I),$ there exist $p,q\in I$ such that $pq\not\in\phi(I)$. Thus, $(a+p)(b+q)\in I-\phi(I)$, since $ab+aq+pb\in\phi(I).$ Hence $s(a+p)\in I$ or $s(b+q)\in\delta(I)$. Consequently, we conclude that $I$ is a $\delta$-$S$-primary ideal of $R$ associated to $s$.\hfill$\blacksquare$\\ \\
Let $R$ be a commutative ring with unity, $I$ a proper ideal of $R$. Then, by Proposition~\ref{2amr-2} and by taking $S=\{1\}$, the following results hold.
\begin{rmk}\textnormal{$\ $\\
(1) If $I$ is a weakly prime ideal of $R$ such that $I$ is not a prime ideal, then $I^2=0$ (it suffices to take $\delta=\delta_0$, $\phi=\phi_0$).\\
(2) If $I$ is a weakly primary ideal of $R$ such that $I$ is not a primary ideal, then $I^2=0$ ( it suffices to take $\delta=\delta_1$, $\phi=\phi_0$).\\
(3) If $I$ is an $n$-almost primary ideal of $R$ such that $I$ is not a primary ideal, then $I^2=I^n$ (it suffices to take $\delta=\delta_1$, $\phi=\phi_n$).}
\end{rmk}
Referring to \cite[Proposition 2(4)]{HM}, \cite[Corollary 2.6]{FBT}, it is clear that if $P$ is an $S$-prime ideal (a weakly $S$-prime ideal) of $R$ for some multiplicative subset $S$ of $R$, then there exists an $s\in S$ such that $s\sqrt{0}\subseteq P$ ($P\subseteq\sqrt{0}$ or $s\sqrt{0}\subseteq P$). Hence, by using Proposition~\ref{2amr-2}, we can deduce easily the following corollary which is a generalization of \cite[Proposition 2(4)]{HM}, \cite[Corollary 2.6]{FBT}.
\begin{cor}\textnormal{
Let $I$ be a $\phi$-$\delta$-$S$-primary ideal of $R$ associated to $s\in S.$ Then $I\subseteq\sqrt{\phi(I)}$ or $s\sqrt{\phi(I)}\subseteq\delta(I).$}
\end{cor}
\textbf{Proof.}$\ $\\
Suppose that $I\not\subseteq\sqrt{\phi(I)}$, then $I^2\not\subseteq\phi(I)$ and hence by Proposition~\ref{2amr-2}, $I$ is a $\delta$-$S$-primary ideal of $R$ associated to $s\in S.$ We show that $s\sqrt{\phi(I)}\subseteq\delta(I).$ Suppose on the contrary that $s\sqrt{\phi(I)}\not\subseteq\delta(I),$ then there exists $y\in\sqrt{\phi(I)}$ such that $sy\not\in\delta(I)$. Let $k$ be the minimal positive integer such that $y^k\in\phi(I)\subseteq I,$ since $I$ is a $\delta$-$S$-primary ideal of $R$ associated to $s$, we have $sy\in\delta(I)$ or $sy^{k-1}\in I$ which implies that $sy^{k-1}\in I$, since $sy\not\in\delta(I).$ Again $sy\not\in\delta(I)$ implies that $s^2y^{k-2}\in I.$ Continuing in this process to get that $s^{k-1}y\in I.$ Since $sy\not\in\delta(I)$, we get $s^k\in I$, a contradiction. Hence $s\sqrt{\phi(I)}\subseteq\delta(I).$\hfill$\blacksquare$
\begin{Theorem}\label{2amr-5}\textnormal{
Let $I$ be a proper ideal of $R$. Then the following statements are equivalent.\\
(1) $I$ is a $\phi$-$\delta$-$S$-primary ideal of $R$ associated to $s\in S.$\\
(2) For each  $a\in R$ such that $a\not\in(\delta(I):s)$ we have either $(I:a)\subseteq(I:s)$ or $(I:a)=(\phi(I):a).$\\
(3) For all $A$ and $B$ ideals of $R$, if $AB\subseteq I$ and $AB\not\subseteq\phi(I)$ then $sA\subseteq I$ or $sB\subseteq\delta(I).$}
\end{Theorem}
\textbf{Proof.}$\ $\\
$(1\rightarrow2)$: Let $a\in R$ such that $a\not\in(\delta(I):s)$, then $sa\not\in\delta(I).$ Suppose that $(I:a)\not=(\phi(I):a).$ We show that $(I:a)\subseteq(I:s).$ Let $r\in(I:a)$. If $r\not\in(\phi(I):a)$, then $ra\in I-\phi(I)$ implies that $sr\in I$, since $sa\not\in\delta(I).$ Suppose $r\in(\phi(I):a).$ Since $(I:a)\not=(\phi(I):a),$ let $r'\in(I:a)$ such that $r'\not\in(\phi(I):a).$ So, $ar'\in I-\phi(I)$ implies that $sr'\in I$, since $sa\not\in\delta(I)$. Thus, $a(r+r')\in I-\phi(I)$ implies that $s(r+r')\in I$, since $sa\not\in\delta(I)$. Because $s(r+r')\in I$ and $sr'\in I$ we get $sr\in I$. Consequently, we conclude that $(I:a)\subseteq(I:s).$\\
$(2\rightarrow1)$: Let $a,b\in R$ such that $ab\in I-\phi(I).$ Suppose that $sa\not\in\delta(I)$, we show that $sb\in I$. Since $b\in (I:a)$ and  $b\not\in (\phi(I):a)$, we get immediately, $(I:a)\subseteq(I:s).$ Thus, $b\in (I:s)$ and $sb\in I$. Accordingly, $I$ is a $\phi$-$\delta$-$S$-primary ideal of $R$ associated to $s$.\\
$(2\rightarrow3)$: Let $A$, $B$ be ideals of $R$ such that $AB\subseteq I.$ Suppose that $sA\not\subseteq I$ and $sB\not\subseteq\delta(I).$ We claim that $AB\subseteq\phi(I).$ Let $b\in B-(\delta(I):s)$, then $(I:b)\subseteq(I:s)$ or $(I:b)=(\phi(I):b).$ Since $A\subseteq(I:b)$ and $A\not\subseteq(I:s)$ we get immediately, $(I:b)=(\phi(I):b).$ So, $Ab\subseteq\phi(I)$. For any $c\in B\cap(\delta(I):s),$ then $b+c\in B-(\delta(I):s)$ implies that $A\subseteq(I:b+c)=(\phi(I):b+c).$ Thus, $A(b+c)\subseteq\phi(I)$ implies that $Ac\subseteq\phi(I),$ since $Ab\subseteq\phi(I).$ Consequently, we conclude that $AB\subseteq\phi(I)$.\\
$(3\rightarrow1)$: Let $a,b\in R$ such that $ab\in I-\phi(I)$. Then $<a><b>\subseteq I$ and $<a><b>\not\subseteq\phi(I)$ implies that $s<a>\subseteq I$ or $s<b>\subseteq\delta(I).$ Thus, $sa\in I$ or $sb\in\delta(I)$. Accordingly, $I$ is a $\phi$-$\delta$-$S$-primary ideal of $R$ associated to $s.$\hfill$\blacksquare$\\ \\
The following result can be proved similar to the previous theorem. Hence, we omit the proof.
\begin{Theorem}\label{2amr-20}\textnormal{
Let $I$ be a proper ideal of $R$. Then the following statements are equivalent.\\
(1) $I$ is a $\phi$-$\delta$-$S$-primary ideal of $R$ associated to $s\in S.$\\
(2) For each  $a\in R$ such that $a\not\in(I:s)$ we have either $(I:a)\subseteq(\delta(I):s)$ or $(I:a)=(\phi(I):a).$\\
(3) For all $A$ and $B$ ideals of $R$, if $AB\subseteq I$ and $AB\not\subseteq\phi(I)$ then $sA\subseteq I$ or $sB\subseteq\delta(I).$}
\end{Theorem}
\begin{Theorem}\textnormal{
Let $P$ be a $\phi$-$\delta$-$S$-primary ideal of $R$ associated to $s\in S.$ If $(\phi(P):a)\subseteq\phi(P:a)$ for each $a\in R-P$, then $(P:a)$ is also $\phi$-$\delta$-$S$-primary ideal of $R$ associated to $s.$}
\end{Theorem}
\textbf{Proof.}$\ $\\
Let $x,y\in R$ such that $xy\in(P:a)-\phi(P:a).$ So, $xya\in P-\phi(P)$ implies that $sxa\in P$ or $sy\in\delta(P).$ Hence, $sx\in (P:a)$ or $sy\in\delta(P)\subseteq(\delta(P):a).$ Thus, $(P:a)$ is a $\phi$-$\delta$-$S$-primary ideal of $R$ associated to $s.$\hfill$\blacksquare$
\begin{cor}\textnormal{
Let $P$ be a $\phi$-$\delta$-$S$-primary ideal of $R$ associated to $s\in S$ and $J$ an ideal in $R$ such that $J\not\subseteq P.$ If $(\phi(P):J)\subseteq\phi(P:J)$ and $(\delta(P):J)\subseteq\delta(P:J)$, then $(P:J)$ is a $\phi$-$\delta$-$S$-primary ideal of $R$ associated to $s.$}
\end{cor}
\textbf{Proof.}$\ $\\
Let $a,b\in R$ such that $ab\in(P:J)-\phi(P:J).$ Then $abJ\subseteq P$ and $abJ\not\subseteq\phi(P)$, since $(\phi(P):J)\subseteq\phi(P:J).$ Thus, $<a><b>J\subseteq P$ and $<a><b>J\not\subseteq\phi(P)$ implies that, by Theorem~\ref{2amr-5}, $s<a>\subseteq P\subseteq(P:J)$ or $s<b>J\subseteq\delta(P).$ So, $sa\in P$ or $sb\in(\delta(P):J)\subseteq\delta(P:J).$ Hence $(P:J)$ is a $\phi$-$\delta$-$S$-primary ideal of $R$ associated to $s.$\hfill$\blacksquare$\\\\
Suppose that $I$ is a $\phi$-$\delta$-$S$-primary ideal of $R$ associated to $s\in S$ such that $\phi\not=\phi_\emptyset$. If $(I:s)=(\delta(I):s)$, then the following result holds.
\begin{prp}\label{2amr-6}\textnormal{
Let $I$ be a $\phi$-$\delta$-$S$-primary ideal of $R$ associated to $s\in S$ such that $(I:s)=(\delta(I):s).$ If $I$ is not a $\delta$-$S$-primary ideal, then $sI\sqrt{\phi(I)}\subseteq\phi(I).$}
\end{prp}
\textbf{Proof.}$\ $\\
Since $I$ is a $\phi$-$\delta$-$S$-primary ideal of $R$ associated to $s\in S$ such that $I$ is not a $\delta$-$S$-primary ideal, by Proposition~\ref{2amr-2}, $I^2\subseteq\phi(I).$ Let $a\in\sqrt{\phi(I)},$ if $a\in(I:s)$ then $sa\in I$ and by Proposition~\ref{2amr-2} we get that $saI\subseteq I^2\subseteq\phi(I).$ Therefore we may assume that $a\not\in(I:s)=(\delta(I):s)$, then by Theorem~\ref{2amr-5}, $(I:a)\subseteq(I:s)=(\delta(I):s)$ or $(I:a)=(\phi(I):a).$ Now, if $(I:a)=(\phi(I):a)$, then $aI\subseteq\phi(I)$ implies that $saI\subseteq\phi(I).$ So we may assume that $(I:a)\subseteq(I:s).$ Let $n\geq1$ be the minimal integer such that $a^n\in\phi(I)$. Then $a^{n-1}\in(I:a)\subseteq(I:s)$ implies that $sa^{n-1}\in I.$ Clearly, $n-1\geq2$, since $sa\not\in I.$ If $sa^{n-1}\not\in\phi(I)$, then $sa^{n-1}=(sa)(a^{n-2})\in I-\phi(I)$ implies that $s^2a\in I$ or $sa^{n-2}\in\delta(I)$. But, if $s^2a\in I,$ then $s^2\in(I:a)\subseteq(I:s)$ implies that $s^3\in I$, a contradiction. So, $sa^{n-2}\in\delta(I)$ and hence $a^{n-2}\in(\delta(I):s)=(I:s)$ implies that $sa^{n-2}\in I-\phi(I)$, since $sa^{n-1}\not\in\phi(I)$. Continuing in this process to get that $sa\in I$ which is a contradiction. Therefore, $sa^{n-1}\in\phi(I)$. Let $j$ be the minimal integer such that $sa^j\in\phi(I)$. Then $j>1$, since $sa\not\in\phi(I)$. Suppose there exists $x\in I$ such that $sax\not\in\phi(I).$ Then $sa(a^{j-1}+x)\in I-\phi(I)$ implies that $s^2a\in I$ or $s(a^{j-1}+x)\in\delta(I).$ Since $s^2a\not\in I$, $s(a^{j-1}+x)\in\delta(I)$ which implies that $a^{j-1}+x\in(\delta(I):s)=(I:s).$ Thus, $sa^{j-1}+sx\in I$ implies that $sa^{j-1}\in I$, since $sx\in I$. Since $j>1$ is the minimal integer such that $sa^j\in\phi(I)$, we get $sa^{j-1}\in I-\phi(I)$. Again continuing in this process to get that $sa\in I$ which is a contradiction. Hence, $sax\in\phi(I)$ for each $x\in I$ and for each $a\in\sqrt{\phi(I)}$. Consequently, we conclude that $sI\sqrt{\phi(I)}\subseteq\phi(I).$\hfill$\blacksquare$
\begin{cor}\textnormal{
Let $I$, $J$ be $\phi$-$\delta$-$S$-primary ideals of $R$ associated to $s\in S$ such that $(I:s)=(\delta(I):s)$ and $(I:s)=(\delta(I):s).$ If $I$, $J$ are not $\delta$-$S$-primary and $\phi(J)\subseteq\phi(I),$ then $sIJ\subseteq\phi(I).$}
\end{cor}
\textbf{Proof.}$\ $\\
Since $I$, $J$ are $\phi$-$\delta$-$S$-primary ideals of $R$ not $\delta$-$S$-primary, by Proposition~\ref{2amr-2}, $I\subseteq\sqrt{\phi(I)}$ and $J\subseteq\sqrt{\phi(J)}\subseteq\sqrt{\phi(I)}.$ Thus, by Proposition~\ref{2amr-6}, $sIJ\subseteq sI\sqrt{\phi(I)}=\phi(I).$\hfill$\blacksquare$
\begin{prp}\textnormal{
Let $I$ be a $\phi$-$\delta$-$S$-primary ideal of $R$ associated to $s\in S$ such that $\phi(J)=\phi(I)$ for each ideal $J\subseteq I.$ If $P$ is an ideal in $R$ such that $P\cap S\not=\emptyset$, then $I\cap P$ and $IP$ are $\phi$-$\delta$-$S$-primary ideals of $R$.}
\end{prp}
\textbf{Proof.}$\ $\\
It is clear that $(P\cap I)\cap S=PI\cap S=\emptyset.$ Pick $t\in P\cap S.$ We show that $I\cap P$ is a $\phi$-$\delta$-$S$-primary ideal of $R$ associated to $ts.$ Let $a,b\in R$ such that $ab\in I\cap P-\phi(I\cap P)$, then  $ab\in I\cap P-\phi(I)\subseteq I-\phi(I).$ Thus, $sa\in I$ or $sb\in\delta(I)$ implies that $tsa\in I\cap P$ or $tsb\in\delta(I)\cap\delta(P)=\delta(I\cap P).$ Consequently, $I\cap P$ is a $\phi$-$\delta$-$S$-primary ideal of $R$ associated to $ts.$ We have the similar proof for $IP.$\hfill$\blacksquare$\\ \\
Let $\phi\not=\phi_\emptyset$ be a reduction function of ideals of $R$ such that $\phi(P)=\phi^2(P)$ for each ideal $P$ of $R$. Then the following result holds.
\begin{prp}\label{2amr-7}\textnormal{
The following statements are equivalent.\\
(1) Every $\phi$-$\delta$-$S$-primary ideal of $R$ is a $\delta$-primary.\\
(2) If $I\in\mathfrak{J}(R)$, then $\phi(I)$ is a $\delta$-primary ideal of $R$ and every $\delta$-$S$-primary ideal of $R$ is a $\delta$-primary.}
\end{prp}
\textbf{Proof.}$\ $\\
($1\rightarrow2$): Let $I\in\mathfrak{J}(R)$. From the definition of $\phi$-$\delta$-$S$-primary ideals of $R$ and $\phi(I)=\phi^2(I)$, we have $\phi(I)$ is a $\phi$-$\delta$-$S$-primary ideal of $R$ and every $\phi$-$\delta$-$S$-primary ideal of $R$ is a $\delta$-primary. Hence, $\phi(I)$ is a $\phi(I)$ is a $\delta$-primary ideal of $R$.\\
($2\rightarrow1$): Let $I$ be a $\phi$-$\delta$-$S$-primary ideal of $R$ associated to $s\in S$. It is enough to show that $I$ is a $\delta$-$S$-primary ideal of $R$ associated to $s$. Let $a,b\in R$ such that $ab\in I$. If $ab\not\in\phi(I)$, then $ab\in I-\phi(I)$ implies that $sa\in I$ or $sb\in\delta(I)$. But, if $ab\in\phi(I)$, then $a\in\phi(I)$ implies $sa\in\phi(I)\subseteq I$ or $b\in\delta(\phi(I))\subseteq\delta(I)$ implies $sb\in\delta(I)$, since $\phi(I)$ is a $\delta$-primary ideal of $R$. Thus, $I$ is a $\delta$-$S$-primary ideal of $R$ associated to $s$ hence, by (2), $I$ is a $\delta$-primary.\hfill$\blacksquare$
\begin{cor} \textnormal{
The following assertions are equivalent:\\
(1) Every weakly $S$-prime ideal of $R$ is a prime ideal.\\
(2) $R$ is a domain and every $S$-prime ideal of $R$ is a prime ideal.}
\end{cor}
\textbf{Proof.}$\ $\\
It suffices to take $\phi=\phi_0$ and $\delta=\delta_0$ in Proposition~\ref{2amr-7}.
\hfill$\blacksquare$
\begin{cor} \textnormal{
The following assertions are equivalent:\\
(1) Every weakly $S$-primary ideal of $R$ is a primary ideal.\\
(2) $R$ is a domain and every $S$-primary ideal of $R$ is a primary ideal.}
\end{cor}
\textbf{Proof.}$\ $\\
It suffices to take $\phi=\phi_0$ and $\delta=\delta_1$ in Proposition~\ref{2amr-7}.
\hfill$\blacksquare$
\begin{rmk}\textnormal{
Let $S_1\subseteq S_2$ be multiplicative subsets of $R$ and $I$ an ideal of $R$ disjoint
with $S_2$. Clearly, if $I$ is a $\phi$-$\delta$-$S_1$-primary ideal of $R$ associated to $s\in S_1$, then $I$ is a $\phi$-$\delta$-$S_2$-primary ideal of $R$ associated to $s\in S_2$. However,
the converse is not true in general (take $\phi=\phi_0$, $\delta=\delta_0$ in \cite[Example 2.3]{FBT}).}
\end{rmk}
\begin{prp}\label{2amr-8}\textnormal{
Let $S_1\subseteq S_2$ be multiplicative subsets of $R$ such that for any $s\in S_2$, there exists $t\in S_2$ with $st\in S_1$. If $I$ is a $\phi$-$\delta$-$S_2$-primary ideal of $R$ associated to $s\in S_2$, then $I$ is a  $\phi$-$\delta$-$S_1$-primary ideal of $R$.}
\end{prp}
\textbf{Proof.}$\ $\\
Let $t\in S_2$ such that $st\in S_1.$ We show that $I$ is a $\phi$-$\delta$-$S_1$-primary ideal of $R$ associated to $st\in S_1.$ Let $a,b\in R$ such that $ab\in I-\phi(I),$ then $sa\in I$ implies that $sta\in I$ or $sb\in\delta(I)$ implies that $stb\in\delta(I).$ Consequently, $I$ is a $\phi$-$\delta$-$S_1$-primary ideal of $R$ associated to $st\in S_1.$\hfill$\blacksquare$\\ \\
 Recall that if $S$ is a multiplicative subset of $R$ with $1\in S,$ then $S^\ast=\{r\in R: \frac{r}{1}\in U(S^{-1}R)\}$ is said to be the saturation of $S.$ One can easily see that $S^\ast$ is a multiplicative subset of $R$ containing $S.$ If $S=S^\ast$, then $S$ is called saturated. Moreover, it is clear that $S^{\ast\ast}=S^\ast.$ (See \cite{RG}).
\begin{prp}\textnormal{
$I$ is a $\phi$-$\delta$-$S$-primary ideal of $R$ if and only if $I$ is a  $\phi$-$\delta$-$S^\ast$-primary ideal of $R$.}
\end{prp}
\textbf{Proof.}$\ $\\
First we show that $S^\ast\cap I=\emptyset$. Let $r\in S^\ast\cap I$, then $\frac{r}{1}$ is a unit in $S^{-1}R$, so there exist $a\in R$, $s\in S$ such that $(\frac{r}{1})(\frac{a}{s})=1.$ Hence, there exists $t\in S$ such that $tra=ts$ which implies that $tra\in I\cap S$, a contradiction. Therefore,  $S^\ast\cap I=\emptyset.$ Since $S\subseteq S^\ast$, $I$ is a $\phi$-$\delta$-$S$-primary ideal of $R$ associated to $s\in S$ implies that $I$ is a $\phi$-$\delta$-$S^\ast$-primary ideal of $R$ associated to $s.$ Conversely, suppose that $I$ is a $\phi$-$\delta$-$S^\ast$-primary ideal of $R$ associated to $s\in S^\ast.$ Let $r\in S^\ast$, then $\frac{r}{1}\in U(S^{-1}R)$ implies that $(\frac{r}{1})(\frac{a}{x})=1$, where $a\in R$, $x\in S.$ Hence, there exists $t\in S$ such that $tra=tx\in S$. Take $r'=ta$, then $r'\in S^\ast$ with $r'r=tx\in S.$ Let $S_1=S$, $S_2=S^\ast,$ then, by Proposition~\ref{2amr-8}, $I$ is a $\phi$-$\delta$-$S$-primary ideal of $R.$\hfill$\blacksquare$\\ \\
Recall that $\delta_S(S^{-1}J)=S^{-1}\delta(J)$ and $\phi_S(S^{-1}J)=S^{-1}\phi(J)$ for each $J\in\mathfrak{J}(R).$ Let $I$ be a proper ideal of $R$ such that $\phi(I:a)=(\phi(I):a)$, $\delta(I:a)=(\delta(I):a)$ for each $a\in R$. Moreover, assume that $\delta(S^{-1}I\cap R)=S^{-1}\delta(I)\cap R.$ Then under the two conditions $\phi(I)=(\phi(I):s)$ for some $s\in S$ and $(\phi(I):t)\subseteq(\phi(I):s)$ for each $t\in S,$ the following result holds.
\begin{Theorem}\textnormal{
Let $I$ be a proper ideal of $R$ such that $I\cap S=\emptyset$. Suppose that $\delta_S(S^{-1}I)\not=S^{-1}R$, if $S^{-1}I\not=S^{-1}R$. Then the following statements are equivalent.\\
(1) $I$ is a $\phi$-$\delta$-$S$-primary ideal of $R$ associated to $s\in S.$\\
(2) $(I:s)$ is a $\phi$-$\delta$-primary ideal of $R.$\\
(3) $S^{-1}I$ is a $\phi_S$-$\delta_S$-primary ideal of $S^{-1}R$ and $(I:t)\subseteq(I:s)$ for each $t\in S.$\\
(4) $S^{-1}I$ is a $\phi_S$-$\delta_S$-primary ideal of $S^{-1}R$ and $S^{-1}I\cap R=(I:s).$}
\end{Theorem}
\textbf{Proof.}$\ $\\
Let $I$ be a proper ideal of $R$ such that $I\cap S=\emptyset$. Then $S^{-1}I\not=S^{-1}R$ implies that $\delta_S(S^{-1}I)\not=S^{-1}R.$ Moreover, it is easy to check that $\delta(I)\cap S=\emptyset.$\\
$(1\rightarrow2)$: Since $I\cap S=\emptyset$, $(I:s)\not=R.$ Let $a,b\in R$ such that $ab\in(I:s)-(\phi(I):s)$, then $sab\in I-\phi(I)$ which implies that $s^2a\in I-\phi(I)$ or $sb\in\delta(I).$ Thus, $sa\in I$, since $s^3\not\in\delta(I)$, or $sb\in\delta(I).$ Hence, $a\in(I:s)$ or $b\in(\delta(I):s).$ Consequently, we conclude that $(I:s)$ is a $\phi$-$\delta$-primary ideal of $R.$\\
$(2\rightarrow3)$: $S^{-1}I\not=S^{-1}R$ since $I\cap S=\emptyset.$ Let $\frac{a}{s_1}, \frac{b}{s_2}\in S^{-1}R$ such that $\frac{a}{s_1}\frac{b}{s_2}\in S^{-1}I-\phi_S(S^{-1}I).$ Then
$\frac{ab}{s_1s_2}=\frac{u}{s_3}\in S^{-1}I-\phi_S(S^{-1}I)$ for some $u\in I.$ So there exists $t\in S$ such that $tabs_3=ts_1s_2u\in I.$ If $ts_1s_2u\in\phi(I)$, then $\frac{u}{s_3}\in \phi_S(S^{-1}I),$ a contradiction. Hence,  $tabs_3\in I-\phi(I)$ which implies that $tabs_3\in (I:s)-(\phi(I):s)$, since $\phi(I)=(\phi(I):s).$ Thus, $a\in(I:s)$ or $tbs_3\in(\delta(I):s).$ Therefore, $sa\in I$ implies that $\frac{a}{s_1}\in S^{-1}I$ or $stbs_3\in\delta(I)$ implies that $\frac{b}{s_2}\in S^{-1}\delta(I).$ So, we conclude that  $S^{-1}I$ is a $\phi_S$-$\delta_S$-primary ideal of $S^{-1}R.$ Let $t\in S$ and let $a\in(I:t)$. If $a\in(\phi(I):t)$, then $a\in(\phi(I):s)\subseteq(I:s).$ Therefore, we may assume that $a\not\in(\phi(I):t).$ So, $ta\in I\subseteq(I:s)$ and $ta\not\in\phi(I)=(\phi(I):s).$ Thus, $ta\in(I:s)-(\phi(I):s)$ which implies that $a\in(I:s)$, since $t\not\in(\delta(I):s).$ Consequently, we conclude that $(I:t)\subseteq(I:s).$\\
$(3\rightarrow4)$: By using part(3) we have $S^{-1}I$ is a $\phi_S$-$\delta_S$-primary ideal of $S^{-1}R.$ Let $s\in S$ such that $(I:t)\subseteq(I:s)$ for each $t\in S.$ Then it is easy to check that $(I:s)\subseteq S^{-1}I\cap R.$ Let $a\in S^{-1}I\cap R$, then there exists $t\in S$ such that $ta\in I$. Thus, $a\in(I:t)\subseteq(I:s)$ implies that $S^{-1}I\cap R\subseteq(I:s).$ Hence we conclude that $S^{-1}I\cap R=(I:s).$\\
$(4\rightarrow1)$: Suppose that $S^{-1}I$ is a $\phi_S$-$\delta_S$-primary ideal of $S^{-1}R$ and $S^{-1}I\cap R=(I:s)$ for some $s\in S$. We show that $I$ is a $\phi$-$\delta$-$S$-primary ideal of $R$ associated to $s$. Let $a,b\in R$ such that $ab\in I-\phi(I).$ Then $\frac{a}{1}\frac{b}{1}\in S^{-1}I-\phi_S(S^{-1}I).$ Thus, $\frac{a}{1}\in S^{-1}I$ or $\frac{b}{1}\in \delta_S(S^{-1}I)=S^{-1}\delta(I),$ since $S^{-1}I$ is a $\phi_S$-$\delta_S$-primary ideal of $S^{-1}R.$ If $\frac{a}{1}\in S^{-1}I$ then there exists $t\in S$ such that $ta\in I$ which implies that $a=\frac{ta}{t}\in S^{-1}I\cap R=(I:s).$ So, $sa\in I.$ Similarly, if $\frac{b}{1}\in S^{-1}\delta(I)$ then there exists $t'\in S$ such that $t'b\in\delta(I)$ which implies that $b=\frac{t'b}{t'}\in S^{-1}\delta(I)\cap R=\delta(S^{-1}I\cap R)=\delta(I:s)=(\delta(I):s).$ So, $sb\in\delta(I).$ Hence we conclude that $I$ is a $\phi$-$\delta$-$S$-primary ideal of $R$ associated to $s$.\hfill$\blacksquare$\\ \\
Let $R$ be a ring, $S\subseteq R$ be a multiplicative subset of $R$. Next we give an example of a proper ideal $P$ of $R$ with $P\cap S=\emptyset$ such that if $(\phi(P):s)\not=\phi(P)$ for some $s\in S$, then $P$ is a $\phi$-$\delta$-$S$-primary ideal of $R$ associated to $s$ but $(P:s)$ is not a $\phi$-$\delta$-primary ideal of $R.$
\begin{Example}\textnormal{
Let $R=\mathbb{Z}[x]$, $P=<4x>=4x\mathbb{Z}[x].$ Let $\phi(P)=p^2=<16x^2>$ and $\delta_1(P)=\sqrt{P}=<2x>.$ Let $S=\{2^k:k\geq0\}$. Then it is easy to check that $P\cap S=\emptyset$ and $P$ is an almost-$\delta_1$-$S$-primary ideal of $R$ associated to $s=2\in S$. Also, it is easy to check that $(P^2:s)=(<16x^2>:2)=<8x^2>\not=P^2.$ Moreover, $(P:s)=(<4x>:2)=<2x>$ is not an almost-$\delta_1$-primary ideal of $R$, since $2x\in<2x>-<4x^2>$ but neither $2\in<2x>=\sqrt{<2x>}$ nor $x\in<2x>=\sqrt{<2x>}.$}
\end{Example}
\section{$(\phi,\delta)-(\psi,\gamma)$-Ring Homomorphisms}
 Following to \cite{ST}, let $X$, $Y$ be commutative rings with unities and let $f:X\rightarrow Y$ be a ring homomorphism. Suppose $\delta$, $\phi$ are expansion and reduction functions of ideals of $X$ and $\gamma$, $\psi$ are expansion and reduction functions of ideals of $Y,$ respectively. Then $f$ is said to be $(\delta,\phi)$-$(\gamma,\psi)$-homomorphism if $\delta(f^{-1}(J))=f^{-1}(\gamma(J))$ and $\phi(f^{-1}(J))=f^{-1}(\psi(J))$ for all $J\in\mathfrak{J}(Y).$
\begin{rmk}\textnormal{$\ $\\
(1) If $f:X\rightarrow Y$ is a nonzero epimorphism and $1$ is the unity of $X$, then $f(1)$ is the unity of $Y$.\\
(2) Suppose $f:X\rightarrow Y$ is a nonzero $(\delta,\phi)$-$(\gamma,\psi)$-epimorphism and let $I$ be a proper ideal of $X$ containing $\ker(f)$. Then it is easy to see that $\gamma(f(I))=f(\delta(I))$ and $\psi(f(I))=f(\phi(I)).$ (\cite[Remark 2.11]{ST})\\
(3) If $S$ is a multiplicative subset of $X$ containing 1, then $f(S)$ is a multiplicative subset of $Y$ containing $f(1)$.}
\end{rmk}
\begin{Theorem}\label{3amr-1}\textnormal{
Let $f:X\rightarrow Y$ be a nonzero $(\delta,\phi)$-$(\gamma,\psi)$-epimorphism. Then the following statements are satisfied.\\
(1) If $J$ is a $\psi$-$\gamma$-$f(S)$-primary ideal of $Y$ associated to $f(s)\in f(S)$, then $f^{-1}(J)$ is a $\phi$-$\delta$-$S$-primary ideal of $X$ associated to $s\in S$.\\
(2) If $I$ is a $\phi$-$\delta$-$S$-primary ideal of $X$ associated to $s\in S$ containing $\ker(f)$ and $f$ is surjective, then $f(I)$ is a $\psi$-$\gamma$-$f(S)$-primary ideal of $Y$ associated to $f(s)\in f(S).$}
\end{Theorem}
\textbf{Proof.}$\ $\\
(1) If $S$ is a multiplicative subset of $X$ with $1\in S$, then $f(S)$ is a multiplicative subset of $Y$ with $1=f(1)\in f(S)$, since $f$ is a nonzero epimorphism. Let $J$ be a $\psi$-$\gamma$-$f(S)$-primary ideal of $Y$ associated to $f(s)\in f(S)$. Choose $a,b\in X$ such that $ab\in f^{-1}(J)-\phi(f^{-1}(J))$. Then we have $f(a)f(b)\in J-\psi(J)$. Since $J$ is a $\psi$-$\gamma$-$f(S)$-primary ideal of $Y$ associated to $f(s)\in f(S)$ we conclude that $f(s)f(a)\in J$ or $f(s)f(b)\in\gamma(J)$, which implies that $sa\in f^{-1}(J)$ or $sb\in f^{-1}(\gamma(J))=\delta(f^{-}(J)).$ Hence $f^{-1}(J)$ is a $\phi$-$\delta$-$S$-primary ideal of $X$ associated to $s$.\\
(2) Let $I$ be a a $\phi$-$\delta$-$S$-primary ideal of $X$ associated to $s$ containing $ker(f)$, then the unity in $Y$ is $f(1)\in f(S)$, since $f$ is a nonzero $(\delta,\phi)$-$(\gamma,\psi)$-epimorphism. Choose $x,y\in Y$ such that $xy\in f(I)-\psi(f(I))$. Since $f$ is onto map, we can choose $a,b\in I$ such that $f(a)=x$, $f(b)=y$. This implies that $f(a)f(b)=f(ab)\in f(I)$. Since $\ker(f)\subseteq I$, we conclude that $ab\in I$. If $ab\in\phi(I)$, then $xy=f(ab)\in f(\phi(I))=\psi(f(I))$, which is a contradiction. So, $ab\in I-\phi(I)$. As $I$ is a $\phi$-$\delta$-$S$-primary ideal of $X$ associated to $s$, we have $sa\in I$ or $sb\in\delta(I)$. Thus, we conclude that $f(s)x\in f(I)$ or $f(s)y\in f(\delta(I))=\gamma(f(I)).$ Therefore, $f(I)$ is a $\psi$-$\gamma$-$f(S)$-primary ideal of $Y$ associated to $f(s).$\hfill$\blacksquare$\\ \\
From the above theorem we obtain the following result.
\begin{Theorem}[correspondence theorem]\textnormal{
Let $f:X\rightarrow Y$ be a nonzero $(\delta,\phi)$-$(\gamma,\psi)$-epimorphism. Then $f$ induces to one-to-one correspondence between the $\phi$-$\delta$-$S$-primary ideals of $X$ associated to $s\in S$ containing $\ker(f)$ and the $\psi$-$\gamma$-$f(S)$-primary ideals of $Y$ associated to $f(s)\in f(S)$ in such a way that if $I$ is a $\phi$-$\delta$-$S$-primary ideal of $X$ associated to $s\in S$ containing $\ker(f)$, then $f(I)$ is the corresponding $\psi$-$\gamma$-$f(S)$-primary ideal of $Y$ associated to $f(s)\in f(S),$ and if $J$ is a $\psi$-$\gamma$-$f(S)$-primary ideal of $Y$ associated to $f(s)\in f(S)$, then $f^{-1}(J)$ is the corresponding $\phi$-$\delta$-$S$-primary ideal of $X$ associated to $s\in S$ containing $\ker(f)$.}\hfill$\blacksquare$
\end{Theorem}
Assume that $\delta$, $\phi$ are expansion and reduction functions of ideals of $R,$ respectively. Let $J$ be a proper ideal of $R$ such that $J=\phi(J)$. Then $\gamma:\mathfrak{J}(R/J)\rightarrow\mathfrak{J}(R/J)$ defined by $\gamma(I/J)=\delta(I)/J$, and $\psi:\mathfrak{J}(R/J)\rightarrow\mathfrak{J}(R/J)$ defined by $\psi(I/J)=\phi(I)/J$ are expansion and reduction functions of ideals of $R/J,$, respectively. Moreover, if $S$ is a multiplicative subset of $R$, then $\bar{S}=S/J$ is a multiplicative subset of $R/J,$ where $S/J=\{\bar{s}=s+J\in R/J\ :\ s\in S\}.$\\ \\
Let $Q$ be a proper ideal of $R$, and let $S$ be a multiplicative subset of $R$. Recall that $Q$ is said to be a weakly $\delta$-$S$-primary ideal of $R$ associated to $s\in S$, if whenever $0\not=ab\in Q$ for some $a,b\in R$ then $sa\in Q$ or $sb\in\delta(Q).$
\begin{Theorem}\label{3amr-2}\textnormal{
Let $\delta$, $\phi$ be expansion and reduction functions of ideals of $R$ and let $J$ be a proper ideal of $R$ such that $J=\phi(J)$. For every $L\in\mathfrak{J}(R)$ let $\gamma:\mathfrak{J}(R/J)\rightarrow\mathfrak{J}(R/J)$ be an expansion function of ideals of $R/J$ defined by $\gamma(L+J/J)=\delta(L+J)/J$ and $\psi:\mathfrak{J}(R/J)\rightarrow\mathfrak{J}(R/J)$ be a reduction function of ideals of $R/J$ defined by $\psi(L+J/J)=\phi(L+J)/J$. Then the followings statements hold.\\
(1) A map $f:R\rightarrow R/J$ defined by $f(r)=r+J$ for every $r\in R$ is a $(\delta,\phi)$-$(\gamma,\psi)$-epimorphism.\\
(2) Let $I$ be a proper ideal of $R$ such that $J\subseteq I$, $S$ a multiplicative subset of $R$. Then $I$ is a $\phi$-$\delta$-$S$-primary ideal of $R$ associated to $s\in S$ if and only if $I/J$ is a $\gamma$-$\psi$-$\bar{S}$-primary ideal of $R/J$ associated to $\bar{s}\in\bar{S}.$\\
(3) Let $I$ be a nonzero proper ideal of $R$ such that $\phi^2(I)=\phi(I).$ Then $I$ is a $\phi$-$\delta$-$S$-primary ideal of $R$ associated to $s\in S$ if and only if $I/\phi(I)$ is a weakly $\gamma$-$\bar{S}$-primary ideal of $R/\phi(I)$ associated to $\bar{s}\in\bar{S}.$}
\end{Theorem}
\textbf{Proof.}$\ $\\
(1) It is easy to see that $f$ is a ring-epimorphism with $\ker(f)=J$. Let $K$ be an ideal in $R/J$, then $K=L+J/J$ for some ideal $L\in\mathfrak{J}(R)$. Therefore, $$f^{-1}(\gamma(K))=f^{-1}(\delta(L+J/J))=\delta(L+J)=\delta(f^{-1}(K)),$$
$$f^{-1}(\psi(K))=f^{-1}(\phi(L+J/J))=\phi(L+J)=\phi(f^{-1}(K)),$$
since $f$ is onto. Thus, $f$ is a $(\delta,\phi)$-$(\gamma,\psi)$-epimorphism.\\
(2) Let $I$ be a proper ideal of $R$ such that $J\subseteq I$, $S$ a multiplicative subset of $R$. Since the map $f$ defined in (1) is a $(\delta,\phi)$-$(\gamma,\psi)$-epimorphism with $\ker(f)=J$ and $f(I)=I/J$. Then, by the correspondence theorem, $I$ is a $\phi$-$\delta$-$S$-primary ideal of $R$ associated to $s\in S$ if and only if $I/J$ is a $\gamma$-$\psi$-$\bar{S}$-primary ideal of $R/J$ associated to $\bar{s}\in\bar{S}.$\\
(3) Let $J=\phi(I)$, then $J=\phi(J)$. Moreover, $f(I)=I/\phi(I)$ and $\psi(I/\phi(I))=\phi(I)/\phi(I)=0\in R/\phi(I).$ Hence, by the correspondence theorem, $I$ is a $\phi$-$\delta$-$S$-primary ideal of $R$ associated to $s\in S$ if and only if $I/\phi(I)$ is a weakly $\gamma$-$\bar{S}$-primary ideal of $R/\phi(I)$ associated to $\bar{s}\in\bar{S}.$\hfill$\blacksquare$
\begin{defn}
Let $I$ be a $\phi$-$\delta$-$S$-primary ideal of $R$ not $\delta$-$S$-primary. Then there exist $a,b\in R$ such that $ab\in\phi(I)$ with $sa\not\in I$ and $sb\not\in\delta(I).$ In this case, $(a,b)$ is called a $\phi$-$\delta$-$S$-twin zero of $I.$
\end{defn}
\begin{Theorem}\textnormal{
Suppose that $I$ is a $\phi$-$\delta$-$S$-primary ideal of $R$ associated to $s\in S.$ If there exist $a,b\in R$ such that $(a,b)$ is a $\phi$-$\delta$-$S$-twin zero of $I.$ Then $\sqrt{I}=\sqrt{\phi(I)}.$}
\end{Theorem}
\textbf{Proof.}$\ $\\
Let $a,b\in R$ such that $(a,b)$ is a $\phi$-$\delta$-$S$-twin zero of $I.$ Then $I$ is not a $\delta$-$S$-primary ideal of $R$. Hence, by Proposition~\ref{2amr-2}, $I^2\subseteq\phi(I)\subseteq I$ implies that $\sqrt{I}=\sqrt{\phi(I)}.$\hfill$\blacksquare$
\begin{lemma}\label{3amr-3}\textnormal{
Let $f:X\rightarrow Y$ be a nonzero $(\delta,\phi)$-$(\gamma,\psi)$-epimorphism and let $I$ a $\phi$-$\delta$-$S$-primary ideal of $X$ associated to $s\in S$ such that $\ker(f)\subseteq I.$ Let $a,b\in X$, then $(a,b)$ is a $\phi$-$\delta$-$S$-twin zero of $I$ if and only if $(f(a),f(b))$ is a $\psi$-$\gamma$-$f(S)$-twin zero of $f(I).$}
\end{lemma}
\textbf{Proof.}$\ $\\
By Theorem~\ref{3amr-1}, $f(I)$ is a $\psi$-$\gamma$-$f(S)$-primary ideal of $Y$ associated to $f(s)\in f(S).$ Let $a,b\in R$ such that $(a,b)$ is a $\phi$-$\delta$-$S$-twin zero of $I.$ Then
$ab\in\phi(I)$ with $sa\not\in I$ and $sb\not\in\delta(I).$ So, $f(a)f(b)=f(ab)\in\psi(f(I))$ with $f(s)f(a)\not\in f(I)$, since $\ker(f)\subseteq I$ and $sa\not\in I.$ Similarly, $f(s)f(b)\not\in \gamma(f(I)).$ Thus, $(f(a),f(b))$ is a $\psi$-$\gamma$-$f(S)$-twin zero of $f(I).$ Conversely, let $a,b\in R$ such that $(f(a),f(b))$ is a $\psi$-$\gamma$-$f(S)$-twin zero of $f(I).$ Then $f(a)f(b)=f(ab)\in\psi(f(I))=f(\phi(I))$ with $f(s)f(a)=f(sa)\not\in f(I)$ and $f(s)f(b)\not\in \gamma(f(I))=f(\delta(I)).$ Thus, $ab\in f^{-1}(\psi(f(I)))=\phi(f^{-1}(f(I)))=\phi(I)$, since $\ker(f)\subseteq I.$ Moreover, $sa\not\in f^{-1}(f(I))=I$ and $sb\not\in f^{-1}(\gamma(f(I)))=\delta(I).$ Consequently, we conclude that $(a,b)$ is a $\phi$-$\delta$-$S$-twin zero of $I.$\hfill$\blacksquare$
\begin{cor}\textnormal{
Let $\delta$, $\phi$ be expansion and reduction functions of ideals of $R$ and let $J$ be a proper ideal of $R$ such that $J=\phi(J)$. For every $L\in\mathfrak{J}(R)$ let $\gamma:\mathfrak{J}(R/J)\rightarrow\mathfrak{J}(R/J)$ be an expansion function of ideals of $R/J$ defined by $\gamma(L+J/J)=\delta(L+J)/J$ and $\psi:\mathfrak{J}(R/J)\rightarrow\mathfrak{J}(R/J)$ be a reduction function of ideals of $R/J$ defined by $\psi(L+J/J)=\phi(L+J)/J$. Let $a,b\in R.$ Then the followings statements hold.\\
(1) $(a,b)$ is a $\phi$-$\delta$-$S$-twin zero of $I$ if and only if $(a+J,b+J)$ is a $\psi$-$\gamma$-$\bar{S}$-twin zero of $I/J.$\\
(2) $(a,b)$ is a $\phi$-$\delta$-$S$-twin zero of $I$ if and only if $(a+\phi(I),b+\phi(I))$ is a $\gamma$-$\bar{S}$-twin zero of $I/\phi(I).$}
\end{cor}
\textbf{Proof.}$\ $\\
(1) It follows from Theorem~\ref{3amr-2}(2) and Lemma~\ref{3amr-3}.\\
(2) It follows from Theorem~\ref{3amr-2}(3) and Lemma~\ref{3amr-3}.
\hfill$\blacksquare$
\begin{Theorem}\textnormal{
Suppose that $I$ is a $\phi$-$\delta$-$S$-primary ideal of $R$ associated to $s\in S.$ If there exist $a,b\in R$ such that $(a,b)$ is a $\phi$-$\delta$-$S$-twin zero of $I.$ Then $aI\subseteq \phi(I)$, $bI\subseteq\phi(I).$ In this case $\sqrt{I}=\sqrt{\phi(I)}.$}
\end{Theorem}
\textbf{Proof.}$\ $\\
Since $(a,b)$ is a $\phi$-$\delta$-$S$-twin zero of $I,$ we have $ab\in\phi(I)$, $sa\not\in I$ and $sb\not\in\delta(I)$. So, $I$ is not a $\delta$-$S$-primary ideal of $R$ associated to $s.$ So, by Proposition~\ref{2amr-2}, $I^2\subseteq\phi(I)\subseteq I$ implies that $\sqrt{I}=\sqrt{\phi(I)}.$ Now, we show that $aI\subseteq \phi(I)$ and $bI\subseteq\phi(I)$ case by case.\\
Since $sb\not\in\delta(I)$, by Theorem~\ref{2amr-5}, $(I:b)\subseteq(I:s)$ or $(I:b)=(\phi(I):b).$ Also, since $a\in(I:b)$ and $a\not\in(I:s)$, we get that $(I:b)=(\phi(I):b).$ Hence we conclude that $bI\subseteq\phi(I)$.\\
Similarly, by Theorem~\ref{2amr-20}, $sa\not\in I$ implies that $(I:a)\subseteq(\delta(I):s)$ or $(I:a)=(\phi(I):a).$ Since $b\in(I:a)$ and $b\not\in(\delta(I):s)$, we get that $(I:a)=(\phi(I):a).$ Hence we conclude that $aI\subseteq\phi(I)$.\hfill$\blacksquare$
\begin{defn}
Suppose $I$ is a $\phi$-$\delta$-$S$-primary ideal of $R$ such that $AB\subseteq I$ and $AB\not\subseteq\phi(I)$, where $A$, $B$ are proper ideals of $R$. Then $I$ is said to be a $\phi$-$\delta$-$S$-free twin zero with respect to $AB$ if $(a,b)$ is not a $\phi$-$\delta$-$S$-twin zero of $I$ for every $a\in A$ and $b\in B$. In particular, $I$ is said to be a $\phi$-$\delta$-$S$-free twin zero, if whenever $AB\subseteq I$ with $AB\not\subseteq\phi(I)$, for some ideals $A$, $B$ of $R$, then $\phi$-$\delta$-$S$-free twin zero with respect to $AB.$
\end{defn}
\begin{Theorem}\label{3amr-4}\textnormal{
Let $I$ be a $\phi$-$\delta$-$S$-primary ideal of $R$ associated to $s\in S.$ Then $I$ is a $\phi$-$\delta$-$S$-free twin zero if and only if for ideals $A$, $B$ of $R$ with $AB\subseteq I$ and $AB\not\subseteq\phi(I),$ either $sA\subseteq I$ or $sB\subseteq\delta(I).$}
\end{Theorem}
\textbf{Proof.}$\ $\\
Suppose that $I$ is a $\phi$-$\delta$-$S$-free twin zero, and let $A$, $B$ be ideals of $R$ such that $AB\subseteq I$ and $AB\not\subseteq\phi(I).$ Then $I$ is a $\phi$-$\delta$-$S$-free twin zero with respect to $AB.$ We show that either $sA\subseteq I$ or $sB\subseteq\delta(I).$ Suppose $sB\not\subseteq\delta(I)$. Then there exists $b\in B$ such that $sb\not\in\delta(I).$ Let $a\in A$, then $(a,b)$ is not a $\phi$-$\delta$-$S$-twin zero of $I.$ If $ab\not\in\phi(I)$, then $ab\in I-\phi(I)$ implies that $sa\in I$, since $sb\not\in\delta(I)$. If $ab\in\phi(I)$, then $sa\in I$, since $(a,b)$ is not a $\phi$-$\delta$-$S$-twin zero of $I$ and $sb\not\in\delta(I).$ Accordingly, we conclude that $sA\subseteq I$. Conversely, suppose that if whenever $A$, $B$ are ideals $R$ with $AB\subseteq I$ and $AB\not\subseteq\phi(I),$ then either $sA\subseteq I$ or $sB\subseteq\delta(I).$ We show that $I$ is a $\phi$-$\delta$-$S$-free twin zero. Let $P$, $Q$ be ideals of $R$  with $PQ\subseteq I$ and $PQ\not\subseteq\phi(I).$ Then, by assumption, either $sP\subseteq I$ or $sQ\subseteq\delta(I).$ Let $p\in P$, $q\in Q.$ If $(p,q)$ is a $\phi$-$\delta$-$S$-twin zero of $I$, then $pq\in\phi(I)$ with $sp\not\in I$ and $sq\not\in\delta(I)$, a contradiction, since $sP\subseteq I$ or $sQ\subseteq\delta(I).$ Thus, $(p,q)$ is not a $\phi$-$\delta$-$S$-twin zero of $I$ for every $p\in P$ and $q\in Q$. Hence we conclude that  $I$ is a $\phi$-$\delta$-$S$-free twin zero with respect to $PQ$ which implies that $I$ is a $\phi$-$\delta$-$S$-free twin zero.\hfill$\blacksquare$
\section{$\phi$-$\delta$-$S$-Primary in direct product of rings}
Let $R_i$ be commutative rings with unity for each $i=1, 2$ and $R =R_1\times R_2$ denote the direct product of rings $R_1$, $R_2$. Also, let $S_1$, $S_2$ be   multiplicative subsets of $R_1$, $R_2$ respectively, then $S=S_1\times S_2$ is a  multiplicative subset of $R$. Suppose that $\phi_i$, $\delta_i$ are reduction and expansion functions of ideals of $R_i$ for each $i=1,2$ respectively. Following to \cite{ST}, we define the following two functions: $$\hat{\delta}(I_1\times I_2)=\delta_1(I_1)\times\delta_2(I_2),$$
$$\hat{\phi}(I_1\times I_2)=\phi_1(I_1)\times\phi_2(I_2).$$
Then it is easy to see that $\hat{\delta}$, $\hat{\phi}$ are expansion and reduction functions of ideals of $R$, respectively.
\begin{Theorem}\label{4amr-11}\textnormal{
Let $R_1$ and $R_2$ be commutative rings with $1\not=0$, $R=R_1\times R_2$ a direct product ring, and $S=S_1\times S_2$ a multiplicative subset of $R$. Suppose that $\delta_i$ is an expansion function of ideals of $R_i$ and $\phi_i$ is a reduction function of ideals of $R_i$ for each $i=1,2$ such that $\phi_2(R_2)\not=R_2$. Then the following statements are equivalent\\
(1) $I_1\times R_2$ is a $\hat{\phi}$-$\hat{\delta}$-$S$-primary ideal of $R$ associated to $(s_1,s_2)\in S$.\\
(2) $I_1$ is a $\delta_1$-$S_1$-primary ideal of $R_1$ associated to $s_1$ and $I_1\times R_2$ is a $\hat{\delta}$-$S$-primary ideal of $R$ associated to $(s_1,s_2)$.}
\end{Theorem}
\textbf{Proof.}$\ $\\
$(1\rightarrow2):$ Suppose that $I_1\times R_2$ is a $\hat{\phi}$-$\hat{\delta}$-$S$-primary ideal of $R$ associated to $(s_1,s_2)\in S$ and let $a,b\in R_1$ such that $ab\in I_1$. Then we have $(a,1)(b,1)=(ab,1)\in I_1\times R_2-\hat{\phi}(I_1\times R_2).$ This implies that $(s_1,s_2)(a,1)\in I_1\times R_2$ or $(s_1,s_2)(b,1)\in\hat{\delta}(I_1\times R_2)$. Hence we conclude that $s_1a\in I_1$ or $s_1b\in\delta_1(I_1)$ and thus, $I_1$ is a $\delta_1$-$S_1$-primary ideal of $R_1$ associated to $s_1$. If $I_1\times R_2$ is not a $\hat{\delta}$-$S$-primary ideal of $R$, then by Proposition(1), we have $(I_1\times R_2)^2\subseteq\hat{\phi}(I_1\times R_2)$ which implies that $R_2=\phi_2(R_2)$ a contradiction. Thus, $I_1\times R_2$ is a $\hat{\delta}$-$S$-primary ideal of $R$ associated to $(s_1,s_2)$.\\
$(2\rightarrow1):$ It is clear, since every $\hat{\delta}$-$S$-primary ideal of $R$ associated to $(s_1,s_2)$ is a $\hat{\phi}$-$\hat{\delta}$-$S$-primary ideal.\hfill$\blacksquare$
\begin{Theorem}\label{4amr-12}\textnormal{
Let $R_1$ and $R_2$ be commutative rings with $1\not=0$, $R=R_1\times R_2$ a direct product ring, and $S=S_1\times S_2$ a multiplicative subset of $R$. Suppose that $\delta_i$ is an expansion function of ideals of $R_i$ and $\phi_i$ is a reduction function of ideals of $R_i$ for each $i=1,2$.
Then the following statements are equivalent\\
(1) $I_1\times R_2$ is a $\hat{\phi}$-$\hat{\delta}$-$S$-primary ideal of $R$ associated to $(s_1,s_2)\in S$ that is not $\hat{\delta}$-$S$-primary.\\
(2) $\hat{\phi}(I_1\times R_2)\not=\emptyset$, $\phi_2(R_2)=R_2$ and $I_1$ is a $\phi_1$-$\delta_1$-$S_1$-primary ideal of $R_1$ associated to $s_1$ that is not $\delta_1$-$S_1$-primary.}
\end{Theorem}
\textbf{Proof.}$\ $\\
$(1\rightarrow2):$ Suppose that $I_1\times R_2$ is a $\hat{\phi}$-$\hat{\delta}$-$S$-primary ideal of $R$ associated to $(s_1,s_2)$ that is not $\hat{\delta}$-$S$-primary, then by Proposition(1), we have $(I_1\times R_2)^2\subseteq\hat{\phi}(I_1\times R_2)$ which implies that $\hat{\phi}(I_1\times R_2)\not=\emptyset.$ If $\phi_2(R_2)\not=R_2$, then by Theorem~\ref{4amr-11}, $I_1\times R_2$ is a $\hat{\delta}$-$S$-primary ideal of $R$ associated to $(s_1,s_2)$ which is a contradiction. Thus, $\phi_2(R_2)=R_2$. Moreover, it is easy to see that $I_1$ is a $\phi_1$-$\delta_1$-$S_1$-primary ideal of $R_1$ associated to $s_1$, since $I_1\times R_2$ is a $\hat{\phi}$-$\hat{\delta}$-$S$-primary ideal of $R$ associated to $(s_1,s_2)$. If $I_1$ is a $\delta_1$-$S_1$-primary ideal of $R_1$ associated to $s_1$, then by the correspondence theorem $I_1\times R_2$ is a $\hat{\delta}$-$S$-primary ideal of $R$ associated to $(s_1,s_2)$ which is a contradiction. Hence $I_1$ is a $\phi_1$-$\delta_1$-$S_1$-primary ideal of $R_1$ associated to $s_1$ that is not $\delta_1$-$S_1$-primary.\\
$(2\rightarrow1):$ We show that $I_1\times R_2$ is a $\hat{\phi}$-$\hat{\delta}$-$S$-primary ideal of $R$ associated to $(s_1,s_2)\in S$. Let $(a,c), (b,d)\in R$ such that $(a,c)(b,d)=(ab,cd)\in I_1\times R_2-\hat{\phi}(I_1\times R_2)$. Then $ab\in I_1-\phi_1(I_1)$, since $\phi_2(R_2)=R_2$. This implies that $s_1a\in I_1$ or $s_1b\in\delta_1(I_1),$ and hence we conclude that $(s_1,s_2)(a,c)=(s_1a,s_2c)\in I_1\times R_2$ or $(s_1,s_2)(b,d)=(s_1b,s_2d)\in \delta_1(I_1)\times\delta_2(R_2)=\hat{\delta}(I_1\times R_2)$. Thus, $I_1\times R_2$ is a $\hat{\phi}$-$\hat{\delta}$-$S$-primary ideal of $R$ associated to $(s_1,s_2)$. Finally, if $I_1\times R_2$ is a $\hat{\delta}$-$S$-primary ideal of $R$ associated to $(s_1,s_2)$, then it is easy to see that $I_1$ is a $\delta_1$-$S_1$-primary ideal of $R_1$ associated to $s_1$ which is a contradiction. Hence $I_1\times R_2$ is not a $\hat{\delta}$-$S$-primary ideal of $R$ associated to $(s_1,s_2)$.\hfill$\blacksquare$\\ \\
Now suppose that for each $i=1,2$, if $I_i\not=\phi_i(I_i)$, then $S_i\cap\phi_i(I_i)=\emptyset$ and if $S_i\cap\delta_i(I_i)\not=\emptyset$, then $S_i\cap I_i=S_i\cap\delta_i(I_i)$. Then we obtain the following result.
\begin{Theorem}\label{4amr-1}\textnormal{
Let $R_1$ and $R_2$ be commutative rings with $1\not=0$, $R=R_1\times R_2$ a direct product ring, and $S=S_1\times S_2$ a multiplicative subset of $R$. Suppose that $\delta_i$ is an expansion function of ideals of $R_i$ and $\phi_i$ is a reduction function of ideals of $R_i$ for each $i=1,2$. Let $I=I_1\times I_2$ be a proper ideal of $R$, for some ideals $I_1\not=\phi_1(I_1)$, $I_2\not=\phi_2(I_2)$ of $R_1$, $R_2$, respectively, such that for every $i\in\{1, 2\}$, if $I_i\not=R_i$, then $\delta_i(I_i)\not=R_i$. Then the following statements are equivalent\\
(1) $I$ is a $\hat{\phi}$-$\hat{\delta}$-$S$-primary ideal of $R$ associated to $(s_1,s_2)\in S$.\\
(2) $I_1=R_1$ and $I_2$ is a $\delta_2$-$S_2$-Primary ideal of $R_2$ associated to $s_2$ or $I_2=R_2$ and $I_1$ is a $\delta_1$-$S_1$-Primary ideal of $R_1$ associated to $s_1$ or $s_2\in I_2\cap S_2$ and $I_1$ is a $\delta_1$-$S_1$-Primary ideal of $R_1$ associated to $s_1$ or $s_1\in I_1\cap S_1$ and $I_2$ is a $\delta_2$-$S_2$-Primary ideal of $R_2$ associated to $s_2$.\\
(3) $I$ is a $\hat{\delta}$-$S$-primary ideal of $R$ associated to $(s_1,s_2)\in S$.}
\end{Theorem}
\textbf{Proof.}$\ $\\
$(1\rightarrow2):$ If $I_1=R_1$, then by Theorem~\ref{4amr-11}, $I_2$ is a $\delta_2$-$S_2$-Primary ideal of $R_2$ associated to $s_2$. Similarly, if $I_2=R_2$, then by Theorem~\ref{4amr-11}, $I_1$ is a $\delta_1$-$S_1$-Primary ideal of $R_1$ associated to $s_1$. Assume that $I_1$, $I_2$ are proper ideals of $R_1$, $R_2$, respectively. Let $a\in I_1$, choose $b\in I_2-\phi_2(I_2)$. Then $(a,1)(1,b)=(a,b)\in I-\hat{\phi}(I)$. As $I$ is a $\hat{\phi}$-$\hat{\delta}$-$S$-primary ideal of $R$ associated to $(s_1,s_2)$, we have $(s_1,s_2)(a,1)=(s_1a,s_2)\in I=I_1\times I_2$ or $(s_1,s_2)(1,b)=(s_1,s_2b)\in \hat{\delta}(I)=\delta_1(I_1)\times\delta_2(I_2)$. Thus, $s_2\in S_2\cap I_2$ or $s_1\in S_1\cap \delta_1(I_1)=S_1\cap I_1$. Assume that $s_2\in S_2\cap I_2$. Since $S\cap I=\emptyset$, we have $S_1\cap I_1=\emptyset$. We show that $I_1$ is a $\delta_1$-$S_1$-Primary ideal of $R_1$ associated to $s_1$. Let $a,b\in R_1$ such that $ab\in I_1$, then we have $(a,s_2)(b,1)\in I-\hat{\phi}(I)$, since $s_2\in S_2\cap I_2$ and $s_2\not\in\phi_2(I_2)$. As $I$ is a $\hat{\phi}$-$\hat{\delta}$-$S$-primary ideal of $R$ associated to $(s_1,s_2)\in S$, we have $(s_1,s_2)(a,s_2)=(s_1a,(s_2)^2)\in I$ or $(s_1,s_2)(b,1)=(s_1b,s_2)\in \hat{\delta}(I)$ which implies that $s_1a\in I_1$ or $s_1b\in\delta_1(I_1)$. Thus, $I_1$ is a $\delta_1$-$S_1$-Primary ideal of $R_1$ associated to $s_1$. Similarly, if we assume that $s_1\in S_1\cap I_1$, then $I_2$ is a $\delta_2$-$S_2$-Primary ideal of $R_2$ associated to $s_2$.\\
$(2\rightarrow3):$ If $I_1=R_1$ and $I_2$ is a $\delta_2$-$S_2$-Primary ideal of $R_2$ associated to $s_2$, then by the correspondence theorem, $I$ is a $\hat{\delta}$-$S$-primary ideal of $R$ associated to $(s_1,s_2)$. Similarly, if $I_2=R_2$ and $I_1$ is a $\delta_1$-$S_1$-Primary ideal of $R_1$ associated to $s_1$, then $I$ is a $\hat{\delta}$-$S$-primary ideal of $R$ associated to $(s_1,s_2)$. Now, suppose that $s_1\in I_1\cap S_1$ and $I_2$ is a $\delta_2$-$S_2$-Primary ideal of $R_2$ associated to $s_2$, we show that $I$ is a $\hat{\delta}$-$S$-primary ideal of $R$ associated to $(s_1,s_2)$. Let $(a,c), (b,d)\in R$ such that $(a,c)(b,d)=(ab,cd)\in I$, then $cd\in I_2$ which implies that $s_2c\in I_2$ or $s_2d\in\delta_2(I_2)$. Since $s_1\in S_1\cap I_1$, we have $(s_1,s_2)(a,c)=(s_1a,s_2c)\in I_1\times I_2$ or $(s_1,s_2)(b,d)=(s_1b,s_2d)\in I_1\times \delta_2(I_2)\subseteq \delta_1(I_1)\times \delta_2(I_2)$. Thus, $I$ is a $\hat{\delta}$-$S$-primary ideal of $R$ associated to $(s_1,s_2)$. Similarly, if we assume that $s_2\in S_2\cap I_2$ and $I_1$ is a $\delta_1$-$S_1$-Primary ideal of $R_1$ associated to $s_1$, then $I$ is a $\hat{\delta}$-$S$-primary ideal of $R$ associated to $(s_1,s_2)$.\\
$(3\rightarrow1):$ Clear.\hfill$\blacksquare$\\
Suppose that for each $i=1,2$, if $I_i\not=\phi_i(I_i)$, then $S_i\cap\phi_i(I_i)=\emptyset$ and if $S_i\cap\delta_i(I_i)\not=\emptyset$, then $S_i\cap I_i=S_i\cap\delta_i(I_i)$. Then we obtain the following result.
\begin{Theorem}\label{4amr-2}\textnormal{
Let $R_1$ and $R_2$ be commutative rings with $1\not=0$, $R=R_1\times R_2$. Let $\delta_1$, $\delta_2$ be expansion functions of ideals of $R_1$, $R_2$, respectively, and let $\phi_1$, $\phi_2$ be reduction functions of ideals of $R_1$, $R_2$, respectively. Let $I=I_1\times I_2$ be a proper ideal of $R$, such that $I\not=\hat{\phi}(I)$ for some ideals $I_1$, $I_2$ of $R_1$, $R_2$, respectively, such that for every $i\in\{1, 2\}$, if $I_i\not=R_i$, then $\delta_i(I_i)\not=R_i$. Then $I$ is a $\hat{\phi}$-$\hat{\delta}$-$S$-primary ideal of $R$ associated to $(s_1,s_2)\in S$ that is not $\hat{\delta}$-$S$-primary if and only if one of the following conditions satisfies\\
(1) $I=I_1\times I_2$, where $\phi_1(I_1)\subsetneqq I_1\subsetneqq R_1$, such that $I_1$ is a $\phi_1$-$\delta_1$-$S_1$-primary ideal of $R_1$ associated to $s_1$ that is not $\delta_1$-$S_1$-primary and $I_2=\phi_2(I_2)$ with $s_2\in S_2\cap\phi_2(I_2)$.\\
(2) $I=I_1\times I_2$, where $\phi_2(I_2)\subsetneqq I_2\subsetneqq R_2$, such that $I_2$ is a $\phi_2$-$\delta_2$-$S_2$-primary ideal of $R_2$ associated to $s_2$ that is not $\delta_2$-$S_2$-primary and $I_1=\phi_1(I_1)$ with $s_1\in S_1\cap\phi_1(I_1).$}
\end{Theorem}
\textbf{Proof.}$\ $\\
Suppose $I$ is a $\hat{\phi}$-$\hat{\delta}$-$S$-primary ideal of $R$ associated to $(s_1,s_2)$ that is not $\hat{\delta}$-$S$-primary. Assume that $I_1\not=\phi_1(I_1)$ and $I_2\not=\phi_2(I_2)$, then by Theorem~\ref{4amr-1}, $I$ is a $\hat{\delta}$-$S$-primary ideal of $R$ associated to $(s_1,s_2)$, a contradiction. Therefore $I_1=\phi_1(I_1)$ or $I_2=\phi_2(I_2)$. Without loss of generality we may assume that $I_2=\phi_2(I_2)$. We show that $s_2\in S_2\cap I_2$ or $s_1\in S_1\cap I_1$. Choose $x\in I_1-\phi_1(I_1)$ then for $b\in I_2$ we have $(x,1)(1,b)=(x,b)\in I-\hat{\phi}(I)$. Since $I$ is a $\hat{\phi}$-$\hat{\delta}$-$S$-primary ideal of $R$ associated to $(s_1,s_2)\in S$, we have $(s_1,s_2)(x,1)\in I$ or $(s_1,s_2)(1,b)\in\hat{\delta}(I)$. Therefore, $(s_1x,s_2)\in I=I_1\times I_2$ or $(s_1,s_2b)\in\hat{\delta}(I)=\delta_1(I_1)\times\delta_2(I_2)$ and hence $s_2\in S_2\cap I_2=S_2\cap\phi_2(I_2)$ or $s_1\in S_1\cap \delta_1(I_1)=S_1\cap I_1$.\\
Case(1) Suppose that $s_2\in S_2\cap I_2=S_2\cap\phi_2(I_2)$. Then $S_1\cap I_1=\emptyset$, since $S\cap I=\emptyset$. Next, we show that $I_1$ is a $\phi_1$-$\delta_1$-$S_1$-primary ideal of $R_1$ associated to $s_1$. Observe that $I_1\not=R_1$. For if $I_1=R_1$, then $S\cap I=(S_1\times S_2)\cap(R_1\times\phi_2(I_2))=S_1\times(S_2\cap\phi_2(I_2))\not=\emptyset$, a contradiction. Thus, $I_1\not=R_1$. Let $a,b\in R_1$ such that $ab\in I_1-\phi_1(I_1)$. Then $(a,1)(b,0)=(ab,0)\in I-\hat{\phi}(I)$ implies that $(s_1,s_2)(a,1)=(s_1a,s_2)\in I_1\times I_2$ or $(s_1,s_2)(b,0)=(s_1b,0)\in \delta_1(I_1)\times \delta_2(I_2)$. So, $s_1a\in I_1$ or $s_1b\in \delta_1(I_1)$. Therefore $I_1$ is a $\phi_1$-$\delta_1$-$S_1$-primary ideal of $R_1$ associated to $s_1$. We show that $I_1$ is not a $\delta_1$-$S_1$-primary ideal of $R_1$ associated to $s_1$. Suppose that $I_1$ is a $\delta_1$-$S_1$-primary ideal of $R_1$ associated to $s_1$ and let $(a,c), (b,d)\in R$ such that $(a,c)(b,d)=(ab,cd)\in I$. Then $ab\in I_1$ implies that $s_1a\in I_1$ or $s_1b\in\delta_1(I_1)$. Since $s_2\in S_2\cap\phi_2(I_2)$ then $(s_1,s_2)(a,c)\in I$ or $(s_1,s_2)(b,d)\in\hat{\delta}(I).$ So, $I$ is a $\hat{\delta}$-$S$-primary ideal of $R$ associated to $(s_1,s_2)$, a contradiction. Thus, $I_1$ is a $\phi_1$-$\delta_1$-$S_1$-primary ideal of $R_1$ associated to $s_1$ that is not $\delta_1$-$S_1$-primary.\\
Case(2) Suppose $s_1\in S_1\cap I_1$. Then $S_2\cap I_2=\emptyset$, since $S\cap I=\emptyset$. We show that $I_2=\phi_2(I_2)$ is a $\delta_2$-$S_2$-primary ideal of $R_2$ associated to $s_2$. Let $a,b\in R_2$ such that $ab\in I_2=\phi_2(I_2)$. Choose $x\in I_1-\phi_1(I_1)$, then $(x,a)((1,b)=(x,ab)\in I-\hat{\phi}(I)$. Since $I$ is a $\hat{\phi}$-$\hat{\delta}$-$S$-primary ideal of $R$ associated to $(s_1,s_2)$, we have $(s_1,s_2)(x,a)=(s_1x,s_2a)\in I=I_1\times I_2$ or $(s_1,s_2)(1,b)=(s_1,s_2b)\in\hat{\delta}(I)=\delta_1(I_1)\times\delta_2(I_2)$ which implies that $s_2a\in I_2$ or $s_2b\in\delta_2(I_2)$. Thus, $I_2=\phi_2(I_2)$ is a $\delta_2$-$S_2$-primary ideal of $R_2$ associated to $s_2$. Now, we show that case(2) can't be happened by proving that  $I$ will be a $\hat{\delta}$-$S$-primary ideal of $R$ associated to $(s_1,s_2)$ which is a contradiction. Let $(a,c), (b,d)\in R$ such that $(a,c)(b,d)=(ab,cd)\in I$, then $cd\in I_2$ implies that $s_2c\in I_2$ or $s_2d\in\delta_2(I_2)$. Since $s_1\in S_1\cap I_1$ we have $(s_1,s_2)(a,c)=(s_1a,s_2c)\in I_1\times I_2=I$ or $(s_1,s_2)(b,d)=(s_1b,s_2d)\in I_1\times\delta_2(I_2)\subseteq\hat{\delta}(I)$. Thus, $I$ is a $\hat{\delta}$-$S$-primary ideal of $R$ associated to $(s_1,s_2)$ which is a contradiction.\\
Conversely, suppose that (1) satisfies. Let $(a,c), (b,d)\in R$ such that $(a,c)(b,d)=(ab,cd)\in I-\hat{\phi}(I)$. Then $ab\in I_1-\phi_1(I_1)$ implies that $s_1a\in I_1$ or $s_1b\in\delta_1(I_1)$. Thus, $(s_1,s_2)(a,c)=(s_1a,s_2c)\in I_1\times I_2=I$ or $(s_1,s_2)(b,d)=(s_1b,s_2d)\in\delta_1(I_1)\times I_2\subseteq\hat{\delta}(I)$. Thus, $I$ is a $\hat{\phi}$-$\hat{\delta}$-$S$-primary ideal of $R$ associated to $(s_1,s_2)$. Finally, we show that $I$ is not a $\hat{\delta}$-$S$-primary ideal of $R$ associated to $(s_1,s_2)$. Suppose that $I$ is a $\hat{\delta}$-$S$-primary ideal of $R$ associated to $(s_1,s_2)$, and let $a,b\in R_1$ such that $ab\in I_1$, then $(a,s_2)(b,1)=(ab,s_2)\in I$. So, $(s_1,s_2)(a,s_2)=(s_1a,(s_2)^2)\in I$ or $(s_1,s_2)(b,1)=(s_1b,s_2)\in\hat{\delta}(I)$ which implies that $s_1a\in I_1$ or $s_1b\in\delta_1(I_1)$. Thus, $I_1$ is a $\delta_1$-$S_1$-primary ideal of $R_1$ associated to $s_1$, a contradiction. Hence $I$ is a $\hat{\phi}$-$\hat{\delta}$-$S$-primary ideal of $R$ associated to $(s_1,s_2)$ that is not $\hat{\delta}$-$S$-primary.\hfill$\blacksquare$
\begin{cor}\textnormal{
Let $R_1$ and $R_2$ be commutative rings with $1\not=0$, $R=R_1\times R_2$. Let $\delta_1$, $\delta_2$ be expansion functions of ideals of $R_1$, $R_2$, respectively. Let $I=I_1\times I_2$ be a proper ideal of $R$ for some ideals $I_1$, $I_2$ of $R_1$, $R_2$, respectively, such that for every $i\in\{1, 2\}$, if $I_i\not=R_i$, then $\delta_i(I_i)\not=R_i$. Then $I$ is a nonzero weakly-$\hat{\delta}$-$S$-primary ideal of $R$ associated to $(s_1,s_2)\in S$ that is not $\hat{\delta}$-$S$-primary if and only if one of the following conditions satisfies\\
(1) $I=I_1\times I_2$, where $I_1$ is a nonzero proper ideal of $R_1$ such that $I_1$ is a weakly-$\delta_1$-$S_1$-primary ideal of $R_1$ associated to $s_1\in S_1$ that is not $\delta_1$-$S_1$-primary and $I_2=0$, $s_2=0$.\\
(2) $I=I_1\times I_2$, where $I_2$ is a nonzero proper ideal of $R_2$ such that $I_2$ is a weakly-$\delta_2$-$S_2$-primary ideal of $R_2$ associated to $s_2\in S_2$ that is not $\delta_2$-$S_2$-primary and $I_1=0$, $s_1=0$.}
\end{cor}
\textbf{Proof.}$\ $\\
In Theorem~\ref{4amr-2}, let $\hat{\phi}(I)=\phi_1(I_1)\times\phi_2(I_2)=(0,0)$ for each proper ideal $I=I_1\times I_2$ of $R$.\hfill$\blacksquare$


\begin{thebibliography}{99}
\bibitem{FBT} \textsl{Almahdi F., Bouba E., Tamekkante M.}, On Weakly $S$-prime ideals of a commutative ring, Analele Stiintifice ale Universitatii Ovidius Constanta, \textbf{29} (2), (2021), 173-186.
\bibitem{AB} \textsl{Anderson D., Bataineh M.} Generalizations of prime ideals, Communications in Algebra 2008, 36 (2): 686-696.
\bibitem{B} \textsl{Badawi A, Fahid B.} On weakly 2-absorbing $\delta$-primary ideals of commutative rings. Georgian Mathematical Journal
2017; 27 (4): 503-516.
\bibitem{RG} \textsl{Gilmer R.}, Multiplicative Ideal Theory, \textbf{12} (M. Dekker, 1972).
\bibitem{HM} \textsl{Hamed A.,  Malek A.}, $S$-prime ideals of a commutative ring, Beitr. Algebra Geom. (2019).
\bibitem{AJ} \textsl{Jaber A.} Properties of $\phi$-$\delta$-primary and 2-absorbing $\delta$-primary ideals of commutative rings, Asian-European Journal of Mathematics (2020), 13 (01): 2050026.
\bibitem{ST2} \textsl{KOC S.}, On weakly 2-prime ideals in commutative rings, Communications in Algebra 2021, 49 (8): 3387-3397.
\bibitem{ST1} \textsl{KOC S., TEKIR U., Yildiz E.}, On weakly 1-absorbing prime ideals, Ricerche di Matematica 2021, 1-16.
\bibitem{ST3} \textsl{Mahdou N., Moutui M., Zahir Y.} Weakly prime ideals issued from an amalgamated algebra, Hacettepe Journal of Mathematics and Statistics \textbf{49} (2020), (3): 1159-1167.
\bibitem{ST} \textsl{YAVUZ S., ONAR S., ERSOY B., TEKIR U., KOC S.}, 2-absorbing
$\phi$-$\delta$-primary ideals, Turkish Journal of Mathematics: \textbf{45}, (2021), 1927-1939.
\bibitem{DZ} \textsl{Zhao D.}, $\delta$-primary Ideals of Commutative Rings,  Kyungpook Mathematical Journal \textbf{41} (2001), 17-22.
\end{thebibliography}
\end{document}